\documentclass[12pt]{article}

\usepackage[a4paper, left=3.5cm,top=3.3cm,right=3.5cm,bottom=4cm]{geometry}

\usepackage{concmath}
\usepackage[T1]{fontenc}

\usepackage{siunitx}
\usepackage{authblk}
\usepackage{psfrag}
\usepackage{amsmath,amssymb,booktabs}

\usepackage[hyphens]{url}
\usepackage{hyperref}
\usepackage{breakurl}

\usepackage{pgfplots}
\definecolor{lime}{RGB}{166,229,117}
\definecolor{goyel}{RGB}{25,25,112}
\definecolor{blor}{RGB}{93,213,249}
\definecolor{blgr}{RGB}{83,154,85}
\definecolor{orred}{RGB}{118,30,176}
\definecolor{drot}{RGB}{160,0,0}
\definecolor{bblau}{RGB}{48,97,172}

\newcommand{\R}{\mathbb R}

\newcommand{\normal}{\boldsymbol n}

\newcommand{\rmd}{\mathrm d}

\newcommand{\coloneqq}{:=}
\newcommand{\la}{\lambda}

\makeatletter
\newcommand*\bigcdot{\mathpalette\bigcdot@{.6}}
\newcommand*\bigcdot@[2]{\mathbin{\vcenter{\hbox{\scalebox{#2}{$\m@th#1\bullet$}}}}}
\makeatother

\newcommand\inner[2]{{#1}\bigcdot {#2}}

\newcommand\abs[1]{\left\vert#1\right\vert}
\newcommand\sabs[1]{\vert#1\vert}
\newcommand\norm[1]{\left\Vert#1\right\Vert}
\newcommand\snorm[1]{\Vert#1\Vert}
\newcommand{\enorm}{\left\|\;\cdot\;\right\|}
\newcommand\set[1]{\left\{#1\right\}}

\newcommand{\dom}{{\boldsymbol D}}
\newcommand{\rand}{{\boldsymbol S}}

\newcommand{\wave}{\mathcal W}
\newcommand{\back}{\mathcal B}
\newcommand{\smooth}{\mathcal K}

\newcommand{\learn}{\mathcal R}
\newcommand{\dal}{\mathcal B}
\newcommand{\unet}{\mathcal N}
\newcommand{\waved}{\mathcal A}
\newcommand{\Dnum}{\mathcal D}

\newcommand{\cset}{C}

\newcommand{\rr}{\boldsymbol r}
\newcommand{\rrs}{\boldsymbol s}

\newcommand{\source}{f}
\newcommand{\pres}{p}
\newcommand{\data}{g}
\newcommand{\irf}{\boldsymbol\phi}
\newcommand{\psf}{\boldsymbol\Phi}

\newcommand{\YY}{\mathbf{G}}
\newcommand{\XX}{\mathbf{F}}
\newcommand{\ZZ}{\boldsymbol{\xi}}
\newcommand{\HH}{\mathbf{H}}
\newcommand{\wdal}{v}
\newcommand{\WW}{\mathbf{W}}
\newcommand{\WWD}{\mathbf {V}}
\newcommand{\WWU}{\mathbf {U}}
\newcommand{\err}{\mathcal E}

\newcommand{\Ns}{M_{\rrs}}
\newcommand{\Nt}{M_t}

\newcommand{\nx}{n}
\newcommand{\Nx}{N_{\rr}}

\newcommand{\Ntrain}{N}
\newcommand{\ntrain}{n}
\newcommand{\ny}{m}

\usepackage{soul}
\usepackage{xcolor}
\colorlet{lred}{red!40}
\colorlet{lgreen}{green!40}
\colorlet{lblue}{blue!40}

\def\mirror at (#1,#2){\filldraw[color=blor](#1,#2-0.05)--(#1+0.15,#2-0.05)--(#1+0.15,#2-0.1)--(#1+0.25,#2)--(#1+0.15,#2+0.1)--(#1+0.15,#2+0.05)--(#1,#2+0.05)--cycle;}
\def\mirrordown at (#1,#2){\filldraw[color=lime](#1-0.02,#2-0.02)--(#1+0.12,#2-0.16)--(#1+0.085,#2-0.195)--(#1+0.21,#2-0.21)--(#1+0.195,#2-0.085)--(#1+0.16,#2-0.12)--(#1+0.02,#2+0.02)--cycle;}
\def\mirrord at (#1,#2){\filldraw[color=lime](#1-0.05,#2)--(#1-0.05,#2-0.15)--(#1-0.1,#2-0.15)--(#1,#2-0.25)--(#1+0.1,#2-0.15)--(#1+0.05,#2-0.15)--(#1+0.05,#2)--cycle;}
\def\mirroru at (#1,#2){\filldraw[color=blgr](#1-0.05,#2)--(#1-0.05,#2+0.15)--(#1-0.1,#2+0.15)--(#1,#2+0.25)--(#1+0.1,#2+0.15)--(#1+0.05,#2+0.15)--(#1+0.05,#2)--cycle;}
\def\mirrortwo at (#1,#2){\filldraw[color=orred](#1,#2-0.05)--(#1+0.15,#2-0.05)--(#1+0.15,#2-0.1)--(#1+0.25,#2)--(#1+0.15,#2+0.1)--(#1+0.15,#2+0.05)--(#1,#2+0.05)--cycle;}

\allowdisplaybreaks

\title{Real-time photoacoustic projection imaging using deep learning}

\author[1]{Johannes Schwab}
\author[1]{Stephan Antholzer}
\author[2]{Robert Nuster}
\author[1,*]{Markus Haltmeier}

\affil[1]{Department of Mathematics, University of Innsbruck, Technikerstra{\ss}e 13, 6020 Innsbruck, Austria}

\affil[2]{Department of Physics, Karl-Franzens-Universitaet Graz,
Universitaetsplatz 5, 8010 Graz,  Austria}

\affil[*]{Corresponding author: markus.haltmeier@uibk.ac.at}

\date{August 30, 2018}

\begin{document}

\maketitle

\begin{abstract}
Photoacoustic tomography (PAT)  is an emerging and  non-invasive hybrid imaging modality for visualizing light absorbing structures in biological tissue.  The recently invented  PAT systems using arrays of 64 parallel integrating line detectors allow capturing  photoacoustic projection images
in fractions of a second. Standard image formation algorithms for this type of setup
suffer from under-sampling due to the sparse detector array,
blurring  due to the finite impulse response of the detection system, and artifacts
due to the limited detection view.  To address these issues, in this paper we develop a new direct
and non-iterative image reconstruction  framework using deep learning.
The proposed DALnet combines the  universal  backprojection (UBP) using  
dynamic aperture length  (DAL) correction  with a deep convolutional neural network (CNN).
Both subnetworks  contain free parameters that are adjusted in the training phase.    
As demonstrated by simulation and experiment, the DALnet is capable of producing 
high-resolution projection images of  3D structures at a frame rate of over 50 
images per second on a standard PC with NVIDIA TITAN Xp GPU. 
The proposed network is shown to outperform state-of-the-art iterative  total variation 
reconstruction algorithms in terms of reconstruction speed as well as in terms of 
various evaluation metrics.       
\end{abstract}

\section{Introduction}
\label{sec:intr}

Photoacoustic tomography (PAT) beneficially combines  the high contrast of pure optical imaging
and the high spatial resolution of pure ultrasound imaging~\cite{Bea11,Wan09b}.
The basic principle of PAT is  as follows. A semitransparent sample (such as  parts of a human patient)
 is illuminated with short pulses of optical radiation. A fraction  of the optical energy is absorbed inside the sample which causes thermal heating, expansion, and a subsequent acoustic pressure wave  depending on the interior light absorbing structures. The acoustic pressure is measured outside of the  sample and used  to reconstruct an image of its interior.
The standard approach for detecting the acoustic waves in PAT is to  use small piezoelectric detector elements arranged on a surface around the object.  In the recent years, PAT systems using
integrating line detectors  have been  invented  and shown to be a high-resolution alternative to the classical approach of using approximate point detectors. Within such systems the pressure signals  are integrated along the  dimension of the line detectors from which  2D projections of the 3D source can be reconstructed \cite{BurHofPalHalSch05,PalNusHalBur07a}. 
By  collecting  2D projections  from different directions, a 3D image of the PA source can be recovered 
with the 2D inverse Radon transform       \cite{BurBauGruHalPal07,haltmeier2009frequency,paltauf2009photoacoustic}.

As the fabrication of an array of parallel line detectors is challenging, initial experiments
for integrating line detectors have been carried out using a single line sensor that is
sequentially moved around the sample in order to collect sufficient  data.
Recently, PAT systems made of arrays with up to $64$ parallel line detectors have been realized~\cite{GraEtAl14,BauPSP15,paltauf2017piezoelectric,Bauer-Marschallinger:17} 
allowing real time PA projection imaging  of 3D structures.
The number of line detectors influences the spatial resolution of the measured  data  and therefore,  the spatial resolution of the final reconstruction.
To increase the spatial resolution, the whole measurement process can be repeated
for different  detector locations. This, however, again increases the measurement time.
If the number of detector locations is far below the Nyquist rate, standard image reconstruction  algorithms (such as filtered backprojection)  yield low quality images containing streak-like under-sampling artifacts.
Consequently, high spatial resolution would require increasing the number of  line sensors.
In order  to address the spatial under-sampling, compressed sensing techniques for PAT have been developed~\cite{arridge2016accelerated,haltmeier2016compressed,haltmeier2017new,sandbichler2015novel,provost2009application}. These methods require
time-consuming iterative  image reconstruction algorithms and, additionally, assume a sparsity prior  that is not always  strictly satisfied~\cite{haltmeier2016compressed}.
Besides under-sampling artefacts, reconstruction with  standard  algorithms additionally suffers from  limited data artifacts as well as blurring due to the  impulse response function  (IRF) of the sensor elements.

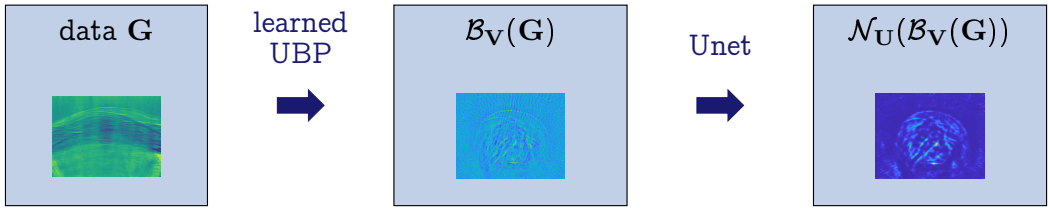
\begin{figure}[thb!]
 \centering
\resizebox{\textwidth}{!}{
\begin{tikzpicture}
\filldraw[color=bblau!30](-0.1,-1.3)--(2.5,-1.3)--(2.5,1.3)--(-0.1,1.3)--cycle;
\draw (-0.1,-1.3)--(2.5,-1.3)--(2.5,1.3)--(-0.1,1.3)--cycle;
\draw (1.2,0.7)node[above]{\small data $\mathbf{G}$};
\node at (1.2,-0.4) {\includegraphics[width=0.1\textwidth]{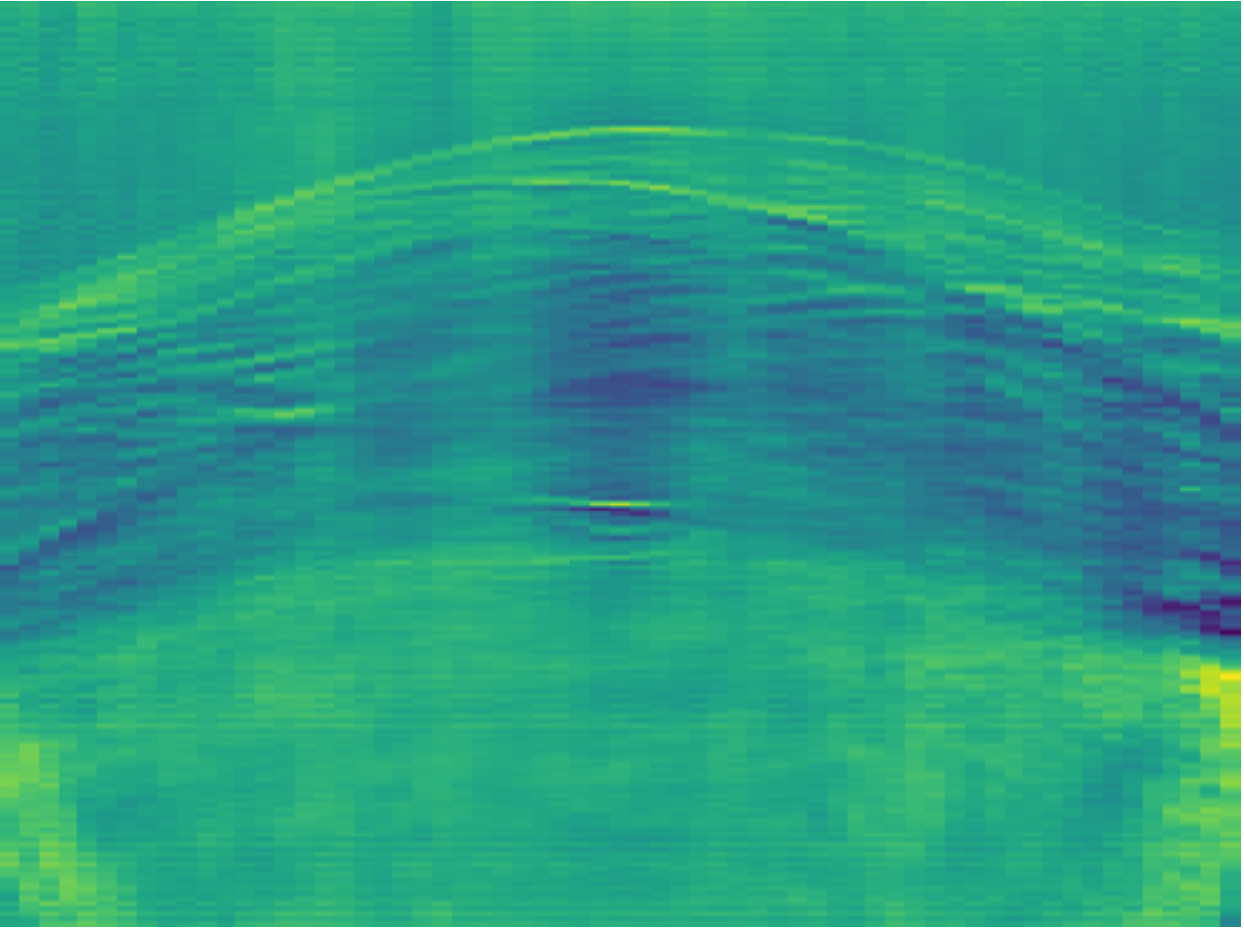}};
\filldraw[color=goyel] (3.4,-0.1)--(3.8,-0.1)--(3.8,-0.2)--(4,0)--(3.8,0.2)--(3.8,0.1)--(3.4,0.1)--cycle;
\draw (3.7,0.8) node[above]{\small \color{goyel}{learned}};
\draw (3.7,0.4) node[above]{\small \color{goyel}{UBP}};
\filldraw[color=bblau!30](4.9,-1.3)--(7.9,-1.3)--(7.9,1.3)--(4.9,1.3)--cycle;
\draw (4.9,-1.3)--(7.9,-1.3)--(7.9,1.3)--(4.9,1.3)--cycle;
\draw (6.4,0.6)node[above]{\small $\dal_{\WWD}(\mathbf{G})$};
\node at (6.4,-0.4) {\includegraphics[width=0.1\textwidth]{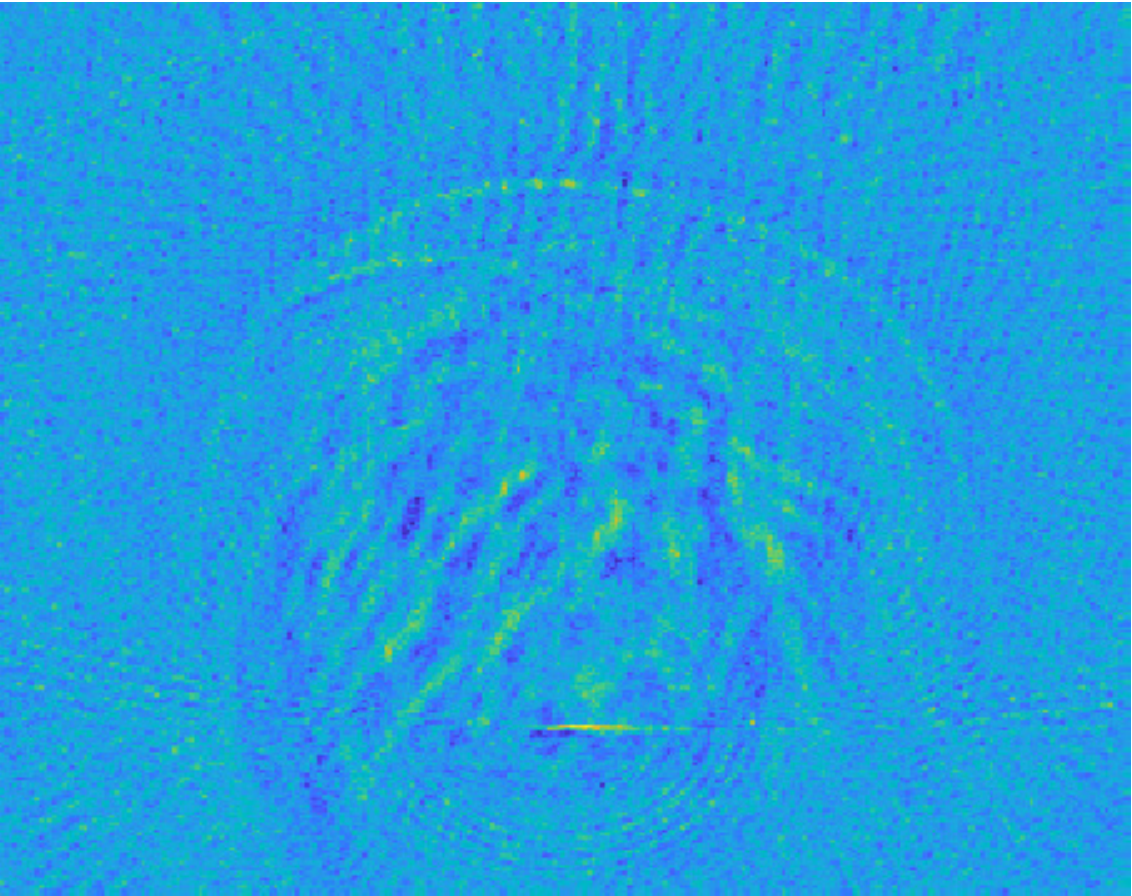}};
\filldraw[color=goyel] (8.8,-0.1)--(9.2,-0.1)--(9.2,-0.2)--(9.4,0)--(9.2,0.2)--(9.2,0.1)--(8.8,0.1)--cycle;
\draw (9.1,0.5) node[above]{\small \color{goyel}{Unet}};
\filldraw[color=bblau!30](10.3,-1.3)--(13.3,-1.3)--(13.3,1.3)--(10.3,1.3)--cycle;
\draw (10.3,-1.3)--(13.3,-1.3)--(13.3,1.3)--(10.3,1.3)--cycle;
\draw (11.8,0.6) node[above]{\small $\unet_{\WWU}(\dal_{\WWD}(\mathbf{G}))$};
\node at (11.8,-0.4) {\includegraphics[width=0.1\textwidth]{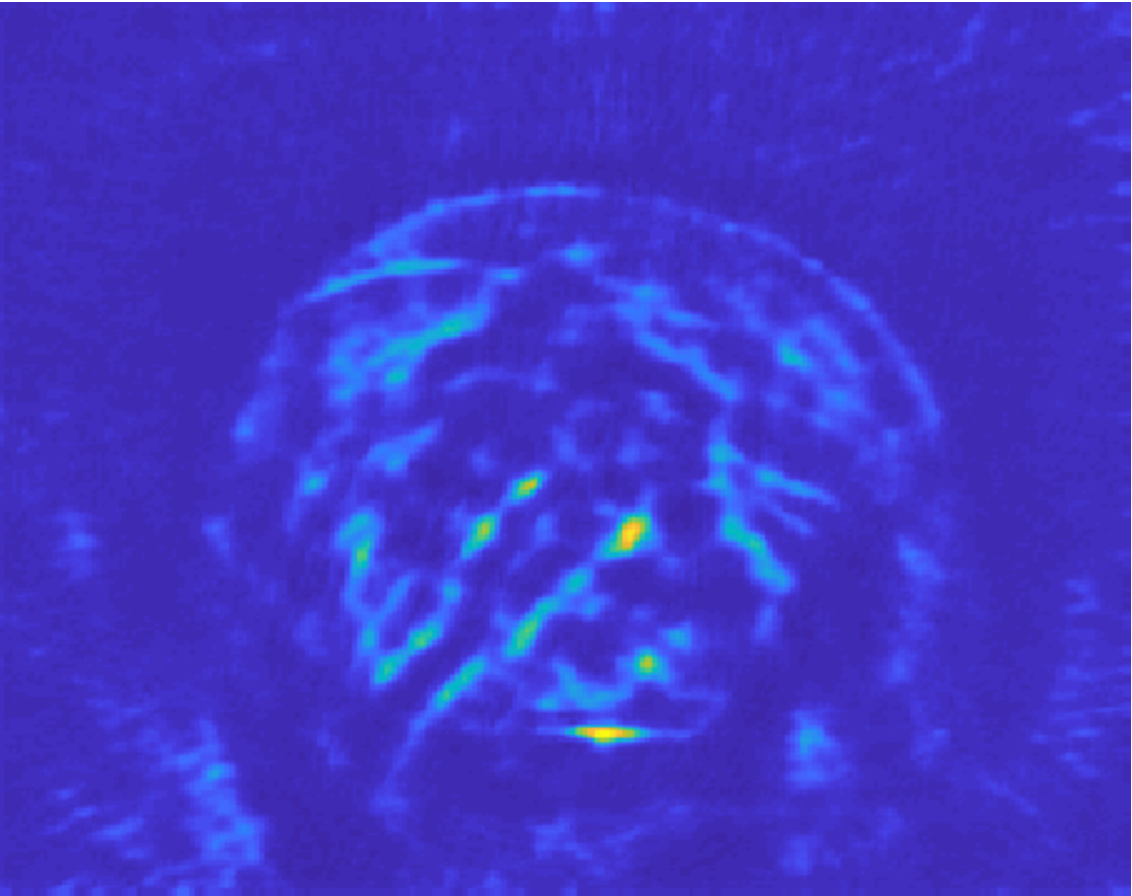}};
\end{tikzpicture}
}
\caption{\textbf{Schematic illustration  of the DALnet.} \label{fig:dalnet} The proposed network consist of a backprojection layer with DAL correction and a subsequent CNN. Both subnetworks contain free parameters  and the DALnet is trained end-to-end.}
\end{figure}

To address the above issues, in  this paper we develop new non-iterative  image reconstruction  algorithms   based on deep learning that are  able to deal with non-ideal IRF, spatial under-sampling and the limited view issue.
We establish a network combining the  universal backprojection (UBP) with dynamic aperture length (DAL) correction \cite{paltauf2007experimental,paltauf2009weight} to address the limited view issue and the IRF   with a convolutional neural network (CNN) to address spatial under-sampling. The resulting DALnet contains  free parameters from the UBP layer and the CNN part which  both are  adjusted end-to-end during the so-called training  phase; see Figure~\ref{fig:dalnet}. 
For sparse sampling PAT,  
image reconstruction with  CNNs has first been proposed in \cite{antholzer2017deep}, where, however, neither the limited view nor the finite bandwidth IRF 
have been  considered.  
Learned filter kernels have be considered  for  x-ray CT in \cite{pelt2014improving}.
CNNs have been previously applied to other medical 
imaging modalities such as (CT) and  magnetic resonance imaging (MRI)
\cite{wang2016perspective,wang2016accelerating,chen2017lowdose,jin2017deep,han2016deep,wurfl2016deep,zhang2016image,rivenson2017deepb}. (For different network architectures used in PAT and other tomographic problems see, for example,  \cite{antholzer2018photoacoustic,allman2018photoacoustic,schwab2018deep,adler2017solving,hauptmann2018model,kelly2017deep,gupta2018cnn,kobler2017variational,waibel2018reconstruction}.)
  However, in all these approaches only the  CNN part contains  trainable  weights, whereas the proposed DALnet  contains trainable weights in the UBP layer and the CNN part.

We train the proposed DALnet  for the PA 64-line detector array system developed in  \cite{paltauf2017piezoelectric} and thereby also account for finite IRF  and the limited view issue.  The algorithm is tested on an experimental  data set consisting of  several PA projection  images of a finger of  one of the  authors.  We demonstrate that real-time high-resolution  PA projection imaging is possible with the trained DALnet. The proposed method is  compared  with the FBP algorithm as well as  state-of-the-art iterative image reconstruction algorithms using total variation (TV) with and without positivity constraint. Our results demonstrate that the 
DALnet  outperforms TV-minimization  in terms of reconstruction speed as well as in various error metrics evaluated on test data not contained in the training set.

\section{Background}
\label{sec:back}

\subsection{Photoacoustic tomography (PAT)}

Suppose some sample of interest is uniformly illuminated with a short optical pulse
which induces an acoustic  pressure wave originating at  light absorbing structures
of the sample.
Assuming constant speed of sound and instantaneous  heating,
the induced acoustic pressure satisfies the free space wave equation 
\begin{equation} \label{eq:wave}
	 \partial_t^2 \pres(\rr, t)  -  v_s^2 \, \Delta_{\rr} \pres(\rr, t)  =
	  \delta'(t) \, \source (\rr) 
	  \,.
	 \end{equation}
Here  $\rr  \in \R^d$ is the spatial location, $t \in \R$ the time, $\Delta_{\rr}$ the spatial Laplacian,
and $v_s$ the speed of sound.  The wave equation  is augmented with  initial condition
$\pres(\rr, t) =0$  for $t < 0$, in which case the  acoustic pressure  is   uniquely  defined
by \eqref{eq:wave}.
The PA source  $\source (\rr)$  is assumed to vanish outside  a bounded
region $\dom \subseteq \R^d$ and
has to be recovered from measurements of the acoustic  pressure outside the sample.

Both cases $d  = 2,3 $  for the  spatial dimension are relevant in PAT.
The case $d=3$ arises in PAT using classical point-wise measurements and  the case $d=2$
for PAT with integrating  line detectors.
In  this work we consider PAT with integrating line detectors and the 2D wave equation.
The PA source  is then a projection image of
the 3D source.  Full  3D  image reconstruction    can be performed by combining projection images from different directions using the 2D inverse Radon transform \cite{BurHofPalHalSch05,PalNusHalBur07a,BurBauGruHalPal07,haltmeier2009frequency,paltauf2009photoacoustic}. However, the PA projection images itself are  already of valuable diagnostic use.

With any PA system, the  acoustic pressure can only be
measured at a finite  number of detector locations, each
having limited bandwidth IRF $\irf$ and recording a finite number of temporal
samples. This yields  to the PAT image reconstruction problem
of reconstructing $\XX$ from data
\begin{equation} \label{eq:data-disc}
     \YY[\ny] =  (\waved(\XX)   \ast_t \irf)[\ny]   + \ZZ[\ny]
      \quad \text{ for }
     \ny \in \set{1, \dots, \Ns} \times \set{1, \dots, \Nt}  \,.
\end{equation}
Here $\XX \in \R^{\Nx \times \Nx}$ is the discretized PA source,
$\waved \colon \R^{\Nx \times \Nx} \to \R^{\Ns \times \Nt}$
the discretization of the solution operator of the wave equation evaluated at the
detector locations, $\irf$ is the discrete IRF and $\ZZ \in \R^{\Ns \times \Nt}$ models the
 noise in the data. The product $\Nx \times \Nx$ is the number of
 reconstruction points, $\Ns$ the number of detector locations and $\Nt$ the number
 of temporal samples. The IRF of the detection system (besides other
 physical factors such as attenuation  and finite optical pulse duration)
 limits the spatial resolution and therefore dictate the number of reconstruction 
 points \cite{XuWan03,HalZan10}. The temporal sampling rate $1/\Nt$ is well above the Nyquist sampling rate
governed by physical  resolution. However, the  number of detector  locations
 is limited in practical applications. Using linear
reconstruction algorithms, the number  $\Ns$ is thereby directly related to the resolution
of the final reconstruction~\cite{haltmeier2016sampling}.

\subsection{Universal backprojection}
\label{sec:ubp}

In the idealized  situation of continuous sampling and full view where the pressure data
$\pres(\rrs,t)$ are known on the whole boundary $\partial  \dom$,
and  assuming an ideal  IRF $\irf(t) = \delta(t)$ (the Dirac delta function) of the
detection system, the inverse problem of recovering the PA source
is uniquely and stably solvable. Several efficient and robust  methods
 including  filtered backprojection  (FBP),
time reversal, or Fourier methods are
available \cite{burgholzer2007exact,Treeby10,HriKucNgu08,FinHalRak07,Hal14,Hal13a,Kun07a,kuchment2008mathematics,rosenthal2013acoustic,xu2005universal} .
A particularly useful  theoretically exact inversion method is the 2D
universal backprojection (UBP) 
\begin{equation}\label{eq:UBP}
          (\back \data)(\rr)
          =
         \frac{1}{\pi}
        \int_{\partial  \dom} \Bigl(
        \inner{\normal_{\rrs}}{(\rr-\rrs)}  \,
         \int_{\abs{\rr-\rrs}}^\infty
        \frac{  \partial_t t^{-1} \data(\rrs,t)}{ \sqrt{t^2-\sabs{\rr-\rrs}^2}}
        \; \rmd t \Bigr)
          \rmd \rrs \,.
\end{equation} 
Here $\data$ is the observed data on the detection curve $\rand$,
$\rr$ is the  reconstruction point,
$\rrs$ the point on the detection curve, and
$\normal_{\rrs}$ the exterior normal to
$\rand \subseteq \partial \dom$ with unit length.

The 2D UBP  has been first derived in \cite{burgholzer2007temporal}
where it is shown to be theoretically exact for linear and circular  geometries, which 
 means that  the PA  source $\source$ can  be 
exactly recovered by the inversion formula $\source = \back \wave \source$, where  
$\wave \source$ denotes the solution of \eqref{eq:wave} restricted to $\partial \dom \times (0, \infty)$.
The 3D version of the UBP has been first derived in \cite{xu2005universal} for planar, spherical and cylindrical geometry. In the recent years, the  UBP has been shown to be theoretically  exact for elliptical and other geometries in arbitrary dimensions~\cite{Hal14,Hal13a,HalPer15b,Kun07a,Nat12}.

\psfrag{X}{$\dom_0$}
\psfrag{Y}{$\dom \setminus \dom_0$}
\psfrag{x}{$\rr$}
\psfrag{a}{$\rrs_+$}
\psfrag{b}{$\rrs_-$}
\psfrag{e}{$\wdal_1$}
\psfrag{f}{$\wdal_2$}
\psfrag{n}{$0$}
\psfrag{1}{$1$}
\psfrag{3}{$\wdal(\rr, \rrs_-)$}
\psfrag{4}{$ \wdal(\rr, \rrs_+) $}
\psfrag{G}{$\rand$}
\begin{figure}[htb!]
    \centering
 \includegraphics[width=0.6\columnwidth]{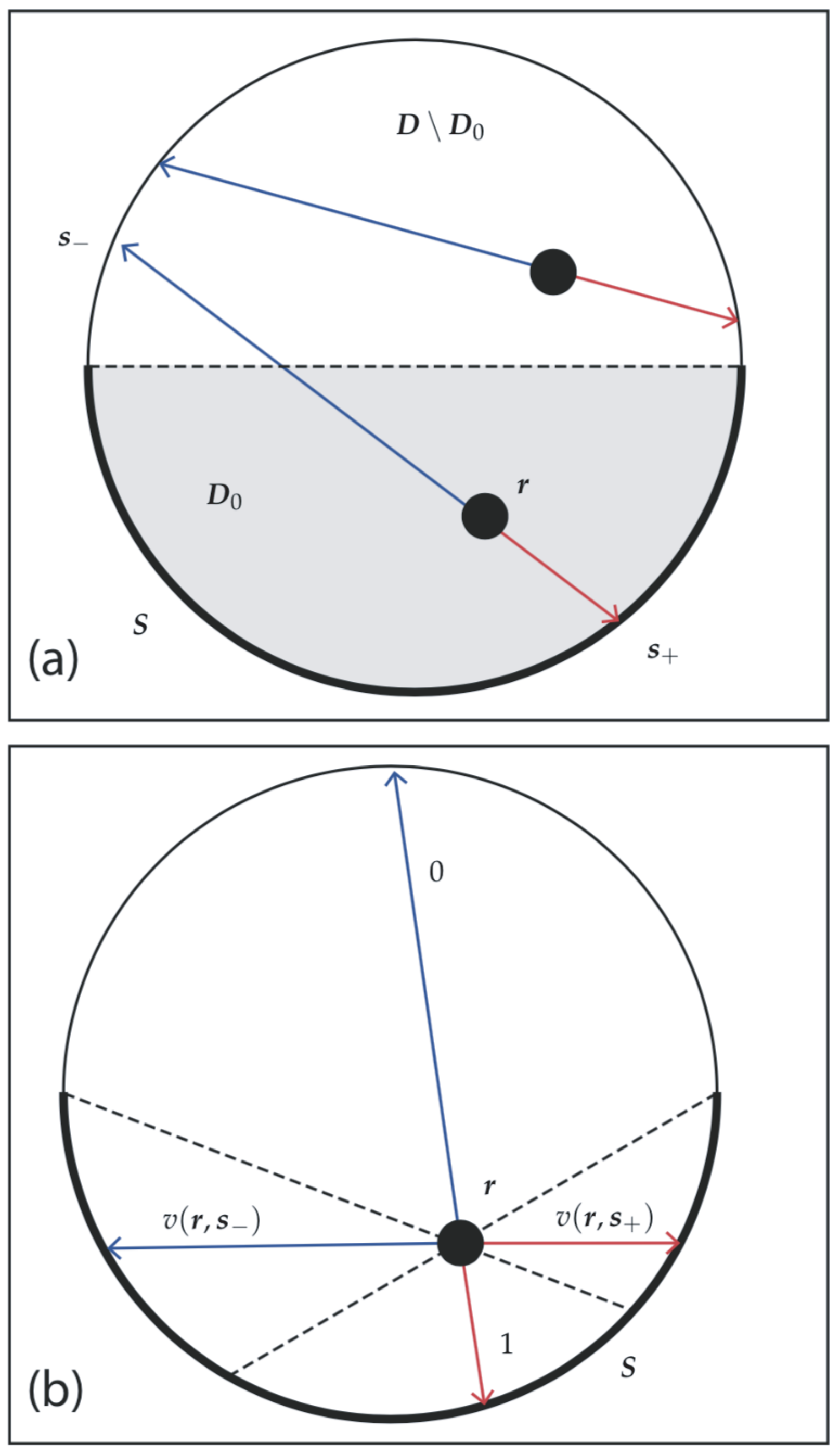}
    \caption{\textbf{Visible zone and UBP with DAL correction.}
    (a) \label{fig:dal}
    Inside the visible zone $\dom_0$ for any pair of rays pointing in
    opposite direction at least one of them
    hits the detection curve.
    (b)  To address the limited view problem, the UBP with DAL correction
    includes a weight with $\wdal(\rr, \rrs_+)  +  \wdal(\rr, \rrs_-) = 1$.
    If $\rrs$ is located outside of  $\rand$, then $\wdal(\rr, \rrs) = 0$.}
\end{figure}

\subsection{Limited view problem and DAL correction}
\label{sec:dal}

In most practical applications, PA data are  available
 only on a proper subset  $\rand  \subsetneq \partial \dom$.
To focus on the main ideas, let us assume that  one  can identify  a  so-called
visible zone  $\dom_0 = \dom_0(\rand) $ defined by the property that  at least one of
any two antipodal rays starting  in the visible zone  hits the detection curve $\rand$;
see Figure~\ref{fig:dal}(a). The reconstruction of $\source$
is known to be stable within  the  visible zone;
however, no direct and theoretically exact inversion method is
available. In such a situation, one can apply iterative methods, where
the forward   and  adjoint problem have  to be solved repeatedly \cite{deanben2012accurate,haltmeier2017iterative,paltauf2002iterative,wang2012investigation}.
This yields  to an increased computation time compared to direct methods.
Typically, the numerical complexity for solving the  forward   and  adjoint
problem is $\mathcal O(\Nx^3)$ which is the same as one application of
the UBP.     For planar geometries Fourier methods  are available  where the
forward and  adjoint  problem (as well as the UBP) can be evaluated with
$\mathcal O(\Nx^2 \log \Nx)$ operation counts.

Another approach, which is computationally less expensive than iterative methods,
is to adapt  the  UBP  to the limited view data.  For full data, any reconstruction point
in the UBP \eqref{eq:UBP} receives information from two antipodal points $\rrs_+$ and  $\rrs_-$.
In the limited data  case, some directions have missing  antipodal points
which yields to blurring of the reconstruction.
To account for this issue, in \cite{paltauf2007experimental} the dynamic aperture length
(DAL)   correction  has been proposed. Using DAL correction, for each  reconstruction point $\rr$,
one considers only  a subpart of the detection curve $\rand$
that has constant view angle $\pi$ from $\rr$ and where exactly one of the two antipodal
detector points is included.  Opposed to simply  extending the  missing  data by zero, the 
DAL correction affects that all singularities are recovered at the correct strength 
\cite{nguyen2014reconstruction,haltmeier2017iterative,SteUhl09}.

In this work, we use a further refinement  of the DAL correction  introduced
 in \cite{paltauf2009weight}, which replaces the UBP formula 
 \eqref{eq:UBP} with the more general weighted UBP formula   
\begin{equation}\label{eq:DAL}
          (\back_{\wdal} \data)(\rr)
          =
        \frac{1}{\pi}
        \int_{\rand}
        \wdal(\rr, \rrs) \Bigl(
        \inner{\normal_{\rrs}}{(\rr-\rrs)}  \,
         \int_{\abs{\rr-\rrs}}^\infty
        \frac{  \partial_t t^{-1} \data(\rrs,t)}{ \sqrt{t^2-\sabs{\rr-\rrs}^2}}
        \; \rmd t \Bigr)
          \rmd \rrs \,.
\end{equation}
The non-negative weight function  $\wdal(\rr, \rrs)$ depending  on the reconstruction  point $\rr$ and the detector location $ \rrs$ is  introduced in such a way that weights for antipodal directions  sum up to one. 
More precisely, the weights satisfy the constraints
$\wdal(\rr, \rrs_+) +   \wdal(\rr, \rrs_-) = 1$  if $\rrs_+, \rrs_-  \in \rand$  are antipodal points  and  $\wdal(\rr, \rrs) = 0$ if $\rrs \in \partial \dom \setminus \rand$; see Figure~\ref{fig:dal}(b).
A  particular deterministic choice for the weight function  $\wdal(\rr, \rrs)$  has been proposed in \cite{paltauf2009weight}.  In this paper compute a weight function   $\wdal(\rr, \rrs)$ 
in  a data driven manner  during the network optimization.

Discretizing \eqref{eq:DAL} yields the approximate left inverse $\dal_{\WWD} \colon  \R^{\Ns \times \Nt}  \to  \R^{\Nx \times \Nx}$ accounting for the limited view issue, where $\WWD$ is a vector including the weights of the DAL. Application of $\dal_{\WWD}$ to  the data in \eqref{eq:data-disc} accounts for the limited view  issue but still yields under-sampling artifacts  as well as blurring due the IRF. The limited view issue is accounted  by   
DAL as it  recovers all singularities at the correct strength assuming sufficient sampling and ideal IRF. 
This implies that even in  case of limited view data, one has 
$\source = \back \wave \source  + \smooth \source $ in  $\dom$, where $\smooth$ 
is at smoothing operator by at least degree one. In the sense of microlocal analysis 
this means that $\back$ is a parametrix for $\wave$ (compare~\cite{nguyen2014reconstruction}). 
Under-sampling and  blurring due to the IRF will be addressed by combining
 $\dal_{\WWD}$ with a convolutional neural network and
deep learning as described in Section~\ref{sec:deep}.

\begin{figure}[htb!]
    \centering
 \includegraphics[width=0.8\columnwidth]{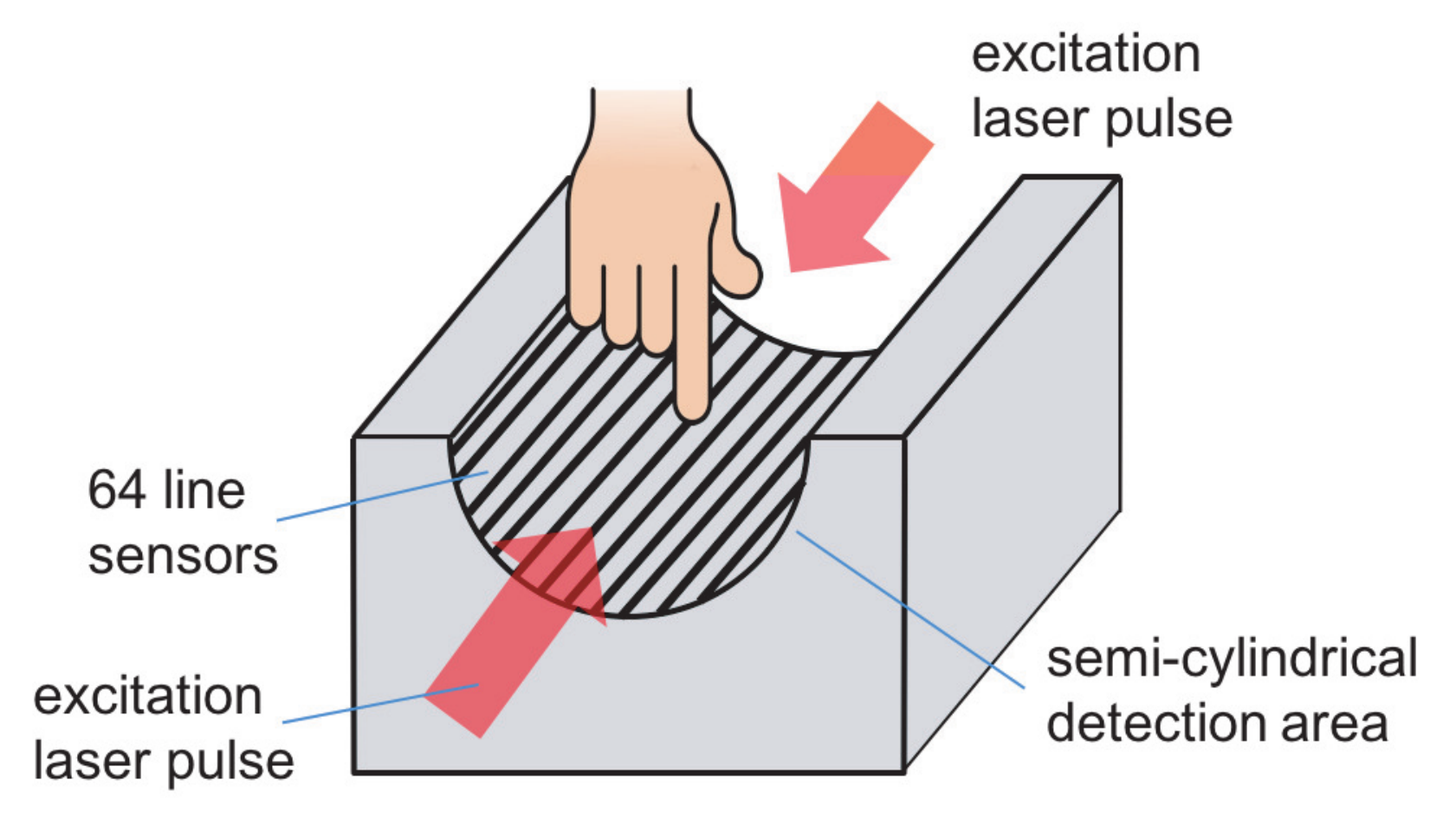}
    \caption{\textbf{Schematic illustration of the tomographic system:}
    The experimental setup contains 64 equispaced piezoelectric line sensors arranged along the detection curve $\rand = \set{\rrs \colon \norm{\rrs}_2 = \SI{50}{mm} \wedge \rrs_2 < 0}$ forming a half-circle, and optical  illumination  is from two sides.    \label{fig:system}}
\end{figure}

\subsection{Tomographic system}
\label{sec:system}

The experimental tomographic system, recently published by Paltauf et.al. \cite{paltauf2017piezoelectric}, was
designed for almost real-time PA projection imaging. In brief, the system is
based on piezoelectric polymer film technology with 64 line shaped sensors homogeneously distributed
on a semi-cylindrical surface with a diameter of \SI{100}{mm}, see Figure~\ref{fig:system}.
The width and length of each sensor is \SI{1.5}{mm} and \SI{150}{mm}, respectively, and the angular increment
between adjacent sensors is \SI{2.8}{\degree}. The design parameters are steered towards
an achievable spatial resolution in the order of $\SI{200}{\micro m}$ to $\SI{250}{\micro m}$ within an imaging area
with radius \SI{20}{mm} centered on the axis of the cylinder. The choice of a half-cylindrical sensor array
arose from considerations regarding the limited-view problem in tomographic reconstruction where all
boundaries of an object can be resolved with maximum resolution if the signals are recorded from at least  \SI{180}{\degree} around the object \cite{paltauf2017piezoelectric}.
The readout of the signals is performed
with a 32-channel data acquisition device, requiring 2:1 multiplexing for each input. Hence, the projection
image acquisition time with the used  \SI{20}{Hz} repetition rate NIR excitation laser system
is \SI{0.1}{seconds} without averaging.

\subsection{Data generation}
\label{sec:data}

We numerically compute data in \eqref{eq:data-disc} at 64 uniformly distributed locations on
the detection curve  $\rand = \{\rrs \colon \norm{\rrs}_2 = \SI{50}{mm} \wedge \rrs_2 < 0\}$.
To generate realistic training data  we also include the experimentally found IRF of the used tomographic  system.
For that purpose we exploit  the convolution relation
$(\waved   \XX) \ast_t \irf = \waved  (\XX  \ast \psf)$ derived in \cite{haltmeier2010spatial}, 
which allows  applying the point  spread function  (PSF)
$\psf$  acting in  the image domain instead of  applying the IRF $\irf$ acting in the
data domain. The explicit relation between $\irf$ and $\psf$ is
computed in~\cite{haltmeier2010spatial}. For the present study, this
explicit relation is  not needed. Instead, we use experimentally measured
data $\YY_{\rm point} \simeq \waved   (\XX_{\rm point}) \ast_t \irf$ corresponding
to a point-like source $\XX_{\rm point}$
to numerically estimate the PSF.  Applying the UBP   and using  the convolution
relation    yields an approximation of the PSF (or convolution kernel) $\psf$.

Having the experimental  PSF at hand,  noisy  data  
\begin{equation} \label{eq:datanoisy}
\YY  = \waved  (\XX  \ast \psf)  + \ZZ
\end{equation}
are computed by  convolving  the PA source $\XX$ with $\psf$, solving the wave equation
\eqref{eq:wave}  and subsequently adding Gaussian white noise $\ZZ$
with standard deviation  \SI{6}{\percent} of the maximum  
$\max [ \waved  (\XX  \ast \psf) ] $ of the exact wave-data.

\section{Image reconstruction}
\label{sec:deep}

\subsection{Proposed DALnet}

Deep learning is a rapidly emerging research field  that  has significantly improved performance of many pattern recognition and machine learning applications~\cite{goodfellow2016deep}.
Deep learning names special  machine learning methods that make use
of deep neural network designs for representing a
nonlinear input to output map together with  optimization procedures
for adjusting the weights of the network during the training phase.
Recently deep learning  algorithms have been developed that
give efficient and highly accurate  tomographic image reconstruction
methods~\cite{antholzer2017deep,chen2017lowdose,han2017deep,jin2017deep,kelly2017deep,kobler2017variational,wang2016accelerating,han2016deep,wang2016perspective,zhang2016image}.
For tomographic image reconstruction, the task of deep learning is
 to find an  image reconstruction function in the form of a deep
 network  $ \learn_\WW $ that maps a measurement data set  $\YY  \in \R^{\Ns \times \Nt}$ to an output image $\XX  \in \R^{\Nx \times \Nx}$.
The weight vector $\WW$ represents a high dimensional   set of free
parameters that can be algorithmically adjusted to the particular reconstruction task.

In our case, the  input images $\YY$  correspond to the noisy measured
PA data \eqref{eq:data-disc} (or \eqref{eq:datanoisy}),
and the output images  $\XX$  are the original  PA sources
evaluated on a discrete grid. The  function $\learn_{\WWU, \WWD}$ is  an 
approximate left   inverse of the PA forward model, that can be well adapted to a relevant
 class of PA sources.
In this article, we propose a network  composed of the  UBP with DAL correction 
to address the limited view issue and  a CNN to address under-sampling
and the finite IRF. Formally, the proposed reconstruction network (DALnet) 
has the form
 \begin{equation} \label{eq:DALnet}
 	\learn_{\WWU, \WWD}    =  \unet_{\WWU} \circ \dal_{\WWD}      \colon \R^{\Ns \times \Nt} \rightarrow\R^{\Nx \times \Nx} \,,
 \end{equation}
where 
\begin{align*}
\dal_{\WWD}.   \colon  \R^{\Ns \times \Nt} \rightarrow\R^{\Nx \times \Nx} \\
\unet_{\WWU}  \colon  \R^{\Nx \times \Nx} \rightarrow \R^{\Nx \times \Nx}
\end{align*}
are  the discretization of the UBP with  DAL correction \eqref{eq:DAL} and
the employed CNN, respectively.

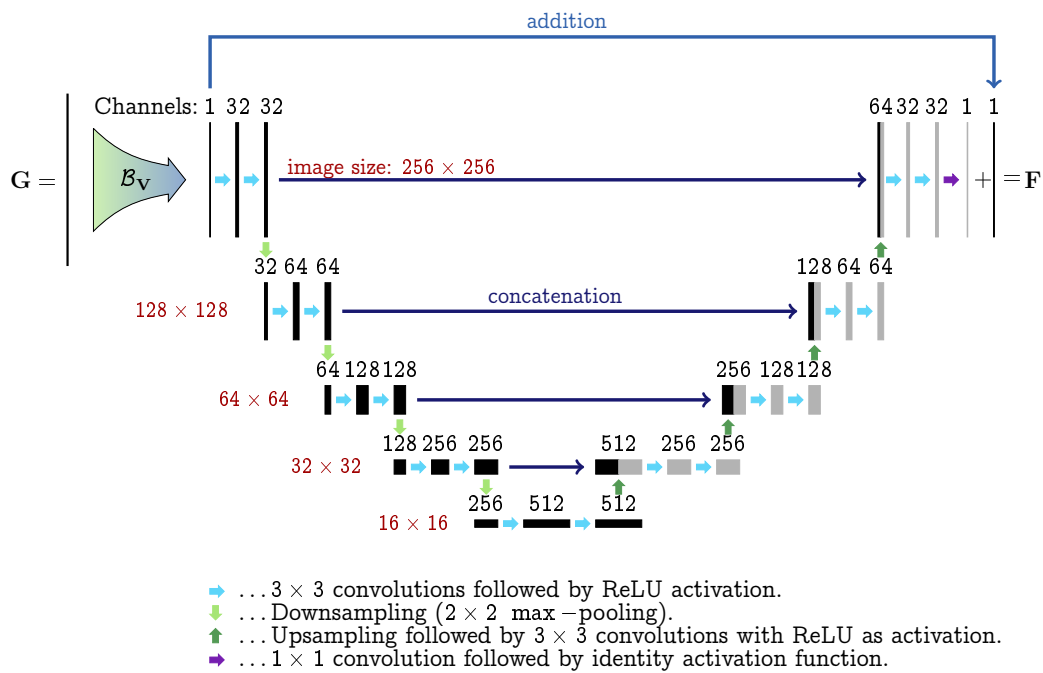
\begin{figure}[thb!]
\centering
\resizebox{\textwidth}{!}{
\begin{tikzpicture}
\draw(0,3)node[left]{};
\draw(0,-7)node[left]{};
\draw(1.2,0.2)node[right]{\small{\color{drot}{image size: $256\times 256$}}};
\filldraw[color=black](0,1)--(0,-1)--(0.01,-1)--(0.01,1)--cycle;
\draw(0.005,1)node[above]{$1$};
\draw(-1.1,1)node[above]{Channels:};
\mirror at (0.1,0)
\filldraw[color=black](0.45,1)--(0.45,-1)--(0.5,-1)--(0.5,1)--cycle;
\draw(0.5,1)node[above]{$32$};
\mirror at (0.6,0)
\filldraw[color=black](0.95,1)--(0.95,-1)--(1,-1)--(1,1)--cycle;
\draw(1.075,1)node[above]{$32$};
\mirrord at (0.975,-1.1)
\draw(0.475,-2.3)node[left]{\small{\color{drot}{$128\times 128$}}};
\filldraw[color=black](0.95,-1.8)--(0.95,-2.8)--(1,-2.8)--(1,-1.8)--cycle;
\draw(0.975,-1.8)node[above]{$32$};
\mirror at (1.1,-2.3)
\filldraw[color=black](1.45,-1.8)--(1.45,-2.8)--(1.55,-2.8)--(1.55,-1.8)--cycle;
\draw(1.5,-1.8)node[above]{$64$};
\mirror at (1.65,-2.3)
\filldraw[color=black](2,-1.8)--(2,-2.8)--(2.1,-2.8)--(2.1,-1.8)--cycle;
\draw(2.05,-1.8)node[above]{$64$};
\mirrord at (2.05,-2.9)
\draw(1.55,-3.85)node[left]{\small{\color{drot}{$64\times 64$}}};
\filldraw[color=black](2,-3.6)--(2,-4.1)--(2.1,-4.1)--(2.1,-3.6)--cycle;
\draw(2.05,-3.6)node[above]{$64$};
\mirror at (2.2,-3.85)
\filldraw[color=black](2.55,-3.6)--(2.55,-4.1)--(2.75,-4.1)--(2.75,-3.6)--cycle;
\draw(2.65,-3.6)node[above]{$128$};
\mirror at (2.85,-3.85)
\filldraw[color=black](3.2,-3.6)--(3.2,-4.1)--(3.4,-4.1)--(3.4,-3.6)--cycle;
\draw(3.3,-3.6)node[above]{$128$};
\mirrord at (3.3,-4.2)
\draw(2.8,-5.025)node[left]{\small{\color{drot}{$32\times 32$}}};
\filldraw[color=black](3.2,-4.9)--(3.2,-5.15)--(3.4,-5.15)--(3.4,-4.9)--cycle;
\draw(3.3,-4.9)node[above]{$128$};
\mirror at (3.5,-5.025)
\filldraw[color=black](3.85,-4.9)--(3.85,-5.15)--(4.15,-5.15)--(4.15,-4.9)--cycle;
\draw(4,-4.9)node[above]{$256$};
\mirror at (4.25,-5.025)
\filldraw[color=black](4.6,-4.9)--(4.6,-5.15)--(5.0,-5.15)--(5.0,-4.9)--cycle;
\draw(4.8,-4.9)node[above]{$256$};
\mirrord at (4.8,-5.25)
\draw(4.3,-6.015)node[left]{\small{\color{drot}{$16\times 16$}}};
\filldraw[color=black](4.6,-5.95)--(4.6,-6.08)--(5,-6.08)--(5,-5.95)--cycle;
\draw(4.8,-5.95)node[above]{$256$};
\mirror at (5.1,-6.015)
\filldraw[color=black](5.45,-5.95)--(5.45,-6.08)--(6.25,-6.08)--(6.25,-5.95)--cycle;
\draw(5.85,-5.95)node[above]{$512$};
\mirror at (6.35,-6.015)
\filldraw[color=black](6.7,-5.95)--(6.7,-6.08)--(7.5,-6.08)--(7.5,-5.95)--cycle;
\draw(7.1,-5.95)node[above]{$512$};
\mirroru at (7.1,-5.5)
\filldraw[color=black](6.7,-4.9)--(6.7,-5.15)--(7.1,-5.15)--(7.1,-4.9)--cycle;
\filldraw[color=black!30](7.1,-4.9)--(7.1,-5.15)--(7.5,-5.15)--(7.5,-4.9)--cycle;
\draw(7.1,-4.9)node[above]{$512$};
\mirror at (7.6,-5.025)
\filldraw[color=black!30](7.95,-4.9)--(7.95,-5.15)--(8.35,-5.15)--(8.35,-4.9)--cycle;
\draw(8.15,-4.9)node[above]{$256$};
\mirror at (8.45,-5.025)
\filldraw[color=black!30](8.8,-4.9)--(8.8,-5.15)--(9.2,-5.15)--(9.2,-4.9)--cycle;
\draw(9,-4.9)node[above]{$256$};
\mirroru at (9,-4.45)
\filldraw[color=black](8.9,-3.6)--(8.9,-4.1)--(9.1,-4.1)--(9.1,-3.6)--cycle;
\filldraw[color=black!30](9.1,-3.6)--(9.1,-4.1)--(9.3,-4.1)--(9.3,-3.6)--cycle;
\draw(9.1,-3.6)node[above]{$256$};
\mirror at (9.4,-3.85)
\filldraw[color=black!30](9.75,-3.6)--(9.75,-4.1)--(9.95,-4.1)--(9.95,-3.6)--cycle;
\draw(9.85,-3.6)node[above]{$128$};
\mirror at (10.05,-3.85)
\filldraw[color=black!30](10.4,-3.6)--(10.4,-4.1)--(10.6,-4.1)--(10.6,-3.6)--cycle;
\draw(10.5,-3.6)node[above]{$128$};
\mirroru at (10.5,-3.15)
\filldraw[color=black](10.4,-1.8)--(10.4,-2.8)--(10.5,-2.8)--(10.5,-1.8)--cycle;
\filldraw[color=black!30](10.5,-1.8)--(10.5,-2.8)--(10.6,-2.8)--(10.6,-1.8)--cycle;
\draw(10.5,-1.8)node[above]{$128$};
\mirror at (10.7,-2.3)
\filldraw[color=black!30](11.05,-1.8)--(11.05,-2.8)--(11.15,-2.8)--(11.15,-1.8)--cycle;
\draw(11.1,-1.8)node[above]{$64$};
\mirror at (11.25,-2.3)
\filldraw[color=black!30](11.6,-1.8)--(11.6,-2.8)--(11.7,-2.8)--(11.7,-1.8)--cycle;
\draw(11.65,-1.8)node[above]{$64$};
\mirroru at (11.65,-1.35)
\filldraw[color=black](11.6,1)--(11.6,-1)--(11.65,-1)--(11.65,1)--cycle;
\filldraw[color=black!30](11.65,1)--(11.65,-1)--(11.7,-1)--(11.7,1)--cycle;
\draw(11.65,1)node[above]{$64$};
\mirror at (11.75,0)
\filldraw[color=black!30](12.1,1)--(12.1,-1)--(12.15,-1)--(12.15,1)--cycle;
\draw(12.125,1)node[above]{$32$};
\mirror at (12.25,0)
\filldraw[color=black!30](12.6,1)--(12.6,-1)--(12.65,-1)--(12.65,1)--cycle;
\draw(12.625,1)node[above]{$32$};
\mirrortwo at (12.75,0)
\filldraw[color=black!30](13.15,1)--(13.15,-1)--(13.16,-1)--(13.16,1)--cycle;
\draw(13.155,1)node[above]{$1$};
\draw(13.12,0)node[right]{$+$};
\filldraw[color=black](13.61,1)--(13.61,-1)--(13.62,-1)--(13.62,1)--cycle;
\draw(13.615,1)node[above]{$1$};
\draw(13.63,0)node[right]{$=$};
\draw(14,0)node[right]{$\mathbf{F}$};
\draw[line width=1.7pt, color=bblau, ->](0.005,1.6)--(0.005,2.5)--(13.615,2.5)--(13.615,1.6);
\draw (6.2,2.5)node[above]{\small \color{bblau}{addition}};
\draw[line width=1.7pt, color=goyel, ->](1.2,0)--(11.4,0);
\draw[line width=1.7pt, color=goyel, ->](2.3,-2.3)--(10.2,-2.3);
\draw (6,-2.3)node[above]{\small \color{goyel}{concatenation}};
\draw[line width=1.7pt, color=goyel, ->](3.6,-3.85)--(8.7,-3.85);
\draw[line width=1.7pt, color=goyel, ->](5.2,-5.025)--(6.5,-5.025);
\mirror at (0,-7.2)
\draw(0.35,-7.2)node[right]{\ldots $3\times 3$ convolutions followed by ReLU activation.};
\mirrord at (0.1,-7.45)
\draw(0.35,-7.6)node[right]{\ldots Downsampling ($2\times 2 \ \max-$pooling).};
\mirroru at (0.1,-8.1)
\draw(0.35,-8)node[right]{\ldots Upsampling followed by $3\times 3$ convolutions with ReLU as activation.};
\mirrortwo at (0,-8.4)
\draw(0.35,-8.4)node[right]{\ldots $1\times 1$ convolution followed by identity activation function.};
\draw(0.35,-8.9)node[right]{};
\filldraw[color=black](-2.48,1.5)--(-2.48,-1.5)--(-2.46,-1.5)--(-2.46,1.5)--cycle;
\draw(-2.5,0)node[left]{$\mathbf{G}=$};
\filldraw[left color=lime!55, right color=bblau!55] (-2.02,0.9) to [out=-90,in=90] (-2.02,-0.9) to [out=45,in=190] (-0.72,-0.25) to (-0.72,-0.38) to (-0.42,0) to (-0.72,0.38) to (-0.72,0.25) to [out=170,in=-45](-2.02,0.9); 
\draw(-0.85,0)node[left]{$\dal_{\WWD}$};
\end{tikzpicture}
}

\caption{\textbf{Architecture of the proposed  DALnet.} The first layer $\dal_{\WWD}$
is a UBP with DAL correction and 
the remaining layers form  a Unet $\unet_{\WWU}$ with residual connection. 
The numbers in red denote the image size,  the numbers above the image stacks 
denote the number of used filters.} \label{fig:unet}
\end{figure}

The CNN part  $\unet_{\WWU} $ can in principle  have  any network structure.
In this work, we use the Unet with residual connection \cite{ronneberger2015unet,han2016deep,jin2017deep,ye2018deep} that we describe in 
more detail in   following subsection. The   weight vector  $ [\WWU, \WWD] \in \R^p$  is  the  vector of adjustable parameters in the DALnet \eqref{eq:DALnet}. The learnable weights $\WWU$ in the 
UBP layer correspond to a discretized  version of the weight  function $\wdal(\rr, \rrs)$. 
We train the  DALnet end-to-end  which means that   
$\WWU, \WWD$  are jointly computed  during the so-called training phase, 
see Subsection~\ref{sec:train}.    
The used  DALnet architecture  is shown in Figure~\ref{fig:unet}.

\subsection{Unet architecture}

The employed CNN part of our reconstruction framework consists of a so-called Unet 
introduced in \cite{ronneberger2015unet}. Originally, the Unet was designed for biomedical image segmentation and improved versions have recently been used for various image reconstruction tasks~\cite{antholzer2017deep,han2016deep,jin2017deep,ye2018deep}. 

The Unet is a fully convolutional network having a encoder-decoder structure with additional skip connections. The encoder part of the Unet consists of four blocks, where in each block the signal is passed through two convolution layers followed by the rectified linear unit (ReLU) activation function defined by $\operatorname{ReLU}(x)\coloneqq \max\{x,0\}$. Due to one max-pooling layer in each block the size of the images decreases while the number of filters increases in every block (see Figure~\ref{fig:unet}). The decoder part of the Unet has the reverse structure, where the downsampling operation (max-pooling) is replaced by a upsampling layer. The skip connections concatenate the interim outputs of the encoder to the upsampled images in the decoder. 
This prevents the network to lose higher frequency parts of the image, which are lost in the down-sampling and additionally helps passing gradient backwardly, finding better local minima~\cite{mao2016image}. Further, there is a  additive skip connection (residual connection), adding the output of the back-projection at the end since the residual images often have simpler structure and training of such residual networks turned out to be more effective~\cite{han2016deep}.

\subsection{Network training}
\label{sec:train}

The  learning aspect  for  the proposed reconstruction network  DALnet defined by \eqref{eq:DALnet}  consists in
adjusting  the parameters $[\WWU, \WWD] \in \R^p$ such that $\learn_{\WWU, \WWD}$
performs well on certain classes of real world data.  For that purpose,
one  chooses pairs of  so-called training data
 \begin{equation} \label{eq:traindata}
 (\YY_\ntrain,\XX_\ntrain) \in \R^{\Nx \times \Nx} \times  \R^{\Ns \times \Nt} 
\quad  \text{ for } \ntrain = 1, \dots,  \Ntrain \,,
 \end{equation} 
where $\XX_\ntrain$ are PA sources  and $\YY_\ntrain$ are the  corresponding noisy  PA data which are computed  according to \eqref{eq:datanoisy}. 
To  measure performance  of the reconstruction network one selects
 an appropriate error function $\err(\WWU, \WWD)$ that quantifies the overall error made by  the network
 $\learn_{\WWU, \WWD}$ on the training data.  During the training  phase the weights are adjusted such that the training error $\err(\WWU, \WWD)$  is minimized.  
 
 Several choices for  the training error are possible. 
 In the present work, we use
 \begin{equation} \label{eq:trainerr}
	\err(\WWU, \WWD) \triangleq
	\frac{1}{2 \Ntrain}
	\sum_{\ntrain=1}^\Ntrain
	 \norm{ \learn_{\WWU, \WWD}( \YY_\ntrain) -  \XX_\ntrain}_2^2 \,,\,
\end{equation}
defined for  all $[\WWU, \WWD] \in \R^p$ where the entries of $\WWD$ are non-negative
and with    $\enorm{}_2$ denoting the  $\ell^2$-norm.  
Standard methods for minimizing  the error function~\eqref{eq:trainerr} are  variants of
stochastic or  incremental gradient  descent algorithms. In this paper we use the stochastic 
(projected) gradient algorithm with momentum \cite{mangasarian1994backpropagation,tseng1998incremental}, where the momentum 
parameter and the step size (learning rate)  are taken as  $0.99$ and  $10^{-3}$, respectively.

\subsection{Constraint  TV minimization}

We  compare the proposed DALnet with the UBP, TV-regularization  and TV-regularization
including a positivity constraint.  For (constraint) TV minimization 
we incorporate the estimated convolution kernel $\psf$ in the forward operator
and construct the output  image as a solution of the optimization
problem
\begin{equation}  \label{eq:TV}
	 \frac{1}{2} \snorm{ \waved  (\XX) \ast  \irf -  \YY}_2^2
	\\ + \la \, \sum_{\nx}  \sqrt{ \abs{ (\Dnum_1 \XX) [\nx] }^2 + \abs{ (\Dnum_2 \XX) [\nx] }^2 }  
	+  I_{\cset} (\XX)
	  \to \min_{\XX}\,.
\end{equation}
Here $\Dnum = [\Dnum_1, \Dnum_2]$ is the discrete gradient operator
and $\la$ the regularization   parameter, and $I_\cset $ denotes the indicator of some convex  set $\cset \subseteq \R^{\Nx \times \Nx}$
defined by $I_\cset(\XX) = 0 $ if  $\XX \in \cset$ and $I_\cset(\XX) =\infty $ else.
In particular, in the case that we take
$\cset = [0, \infty)^{\Nx \times \Nx}$  it guarantees non-negativity. For  
$\cset = \R^{\Nx \times \Nx}$,  \eqref{eq:TV} reduced to standart TV-regularization.    

The discrete  TV problem~\eqref{eq:TV} can be minimized  by various methods.
In this work, we use the minimization algorithm  of \cite{sidky2012convex}, which is
a special  instance of the  Chambolle-Pock algorithm \cite{chambolle2011first}.
The algorithm of \cite{sidky2012convex} for TV minimization has been previously applied 
to  PAT in \cite{boink2018framework,nguyen2018reconstruction}. However, for  
PA projection imaging including the detector PSF it is applied  for the  
first time in the present paper.

\section{Numerical and experimental results}
\label{sec:results}

In this  section, we present reconstruction results  for  UBP, 
TV-minimization with and without positivity constraint and the proposed 
DALnet.  Results are presented for simulated as well as experimental data.

\subsection{Training and evaluation data}
\label{sec:data}

Adequate training and evaluation of the  network is important in order to
perform well on experimental data. For that purpose we
generate 200 projection images which are  computed from a
3D lung blood vessel data set of size  $512 \times 512 \times 512$  by rotating it
along a single axis. From these projection
images, we extracted $3200$ rotated patches $\XX_\ntrain$ of size
$256 \times 256$ for training and evaluation. The corresponding data
are generated by solving the PAT forward problem   according to \eqref{eq:datanoisy} 
corresponding to the 64-line array illustrated in Figure~\ref{fig:system}.

The  generated data  pairs are split in $\Ntrain = 3000$ data pairs for training  and 
200 
data pairs for validation. The DALnet is trained  by minimizing  \eqref{eq:trainerr} as
described above. All training and evaluation patches  are
normalized to have a maximal intensity value of one and cover the imaging
region $[\SI{-12.5}{mm}, \SI{12.5}{mm}] \times [\SI{-20}{mm}, \SI{5}{mm}]$ 
which is  mostly contained inside the detection curve 
$\rand = \set{\rrs \colon \norm{\rrs}_2 = \SI{50}{mm} \wedge \rrs_2 < 0}$.

\begin{figure}[htb!]
    \centering
    \includegraphics[width=0.3\columnwidth, height=0.25\columnwidth]{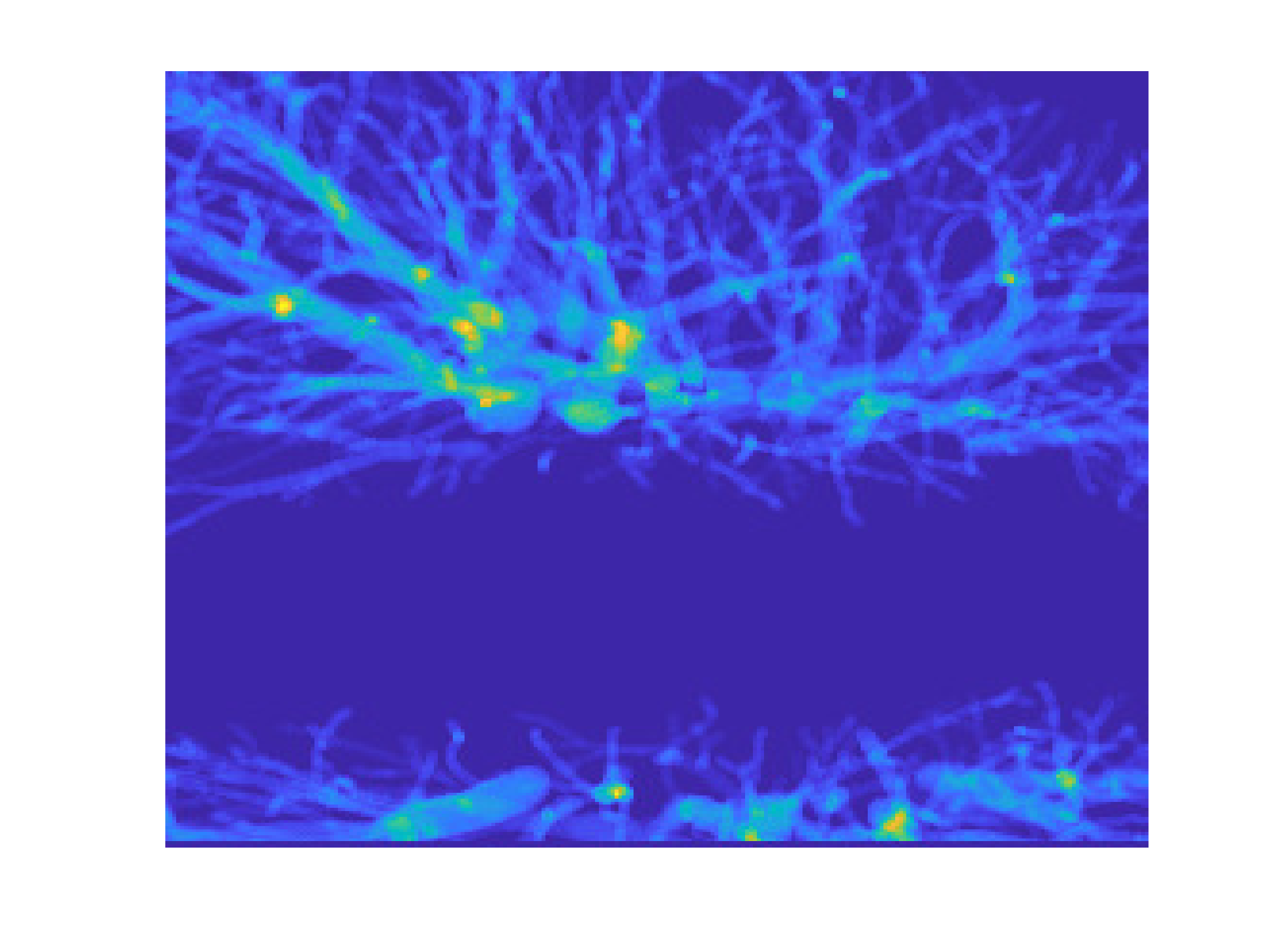}
    \includegraphics[width=0.3\columnwidth, height=0.25\columnwidth]{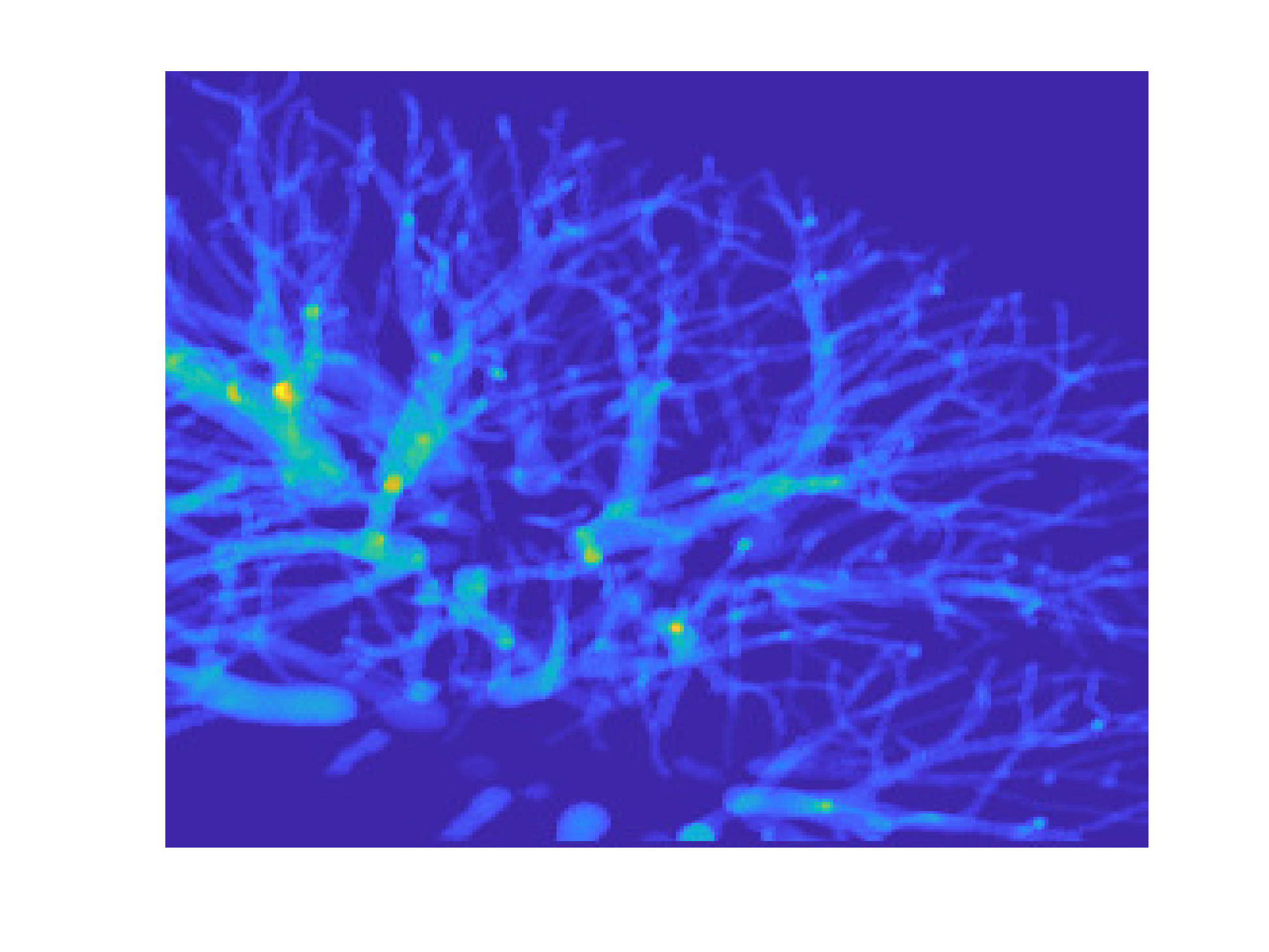}\\
    \includegraphics[width=0.3\columnwidth, height=0.25\columnwidth]{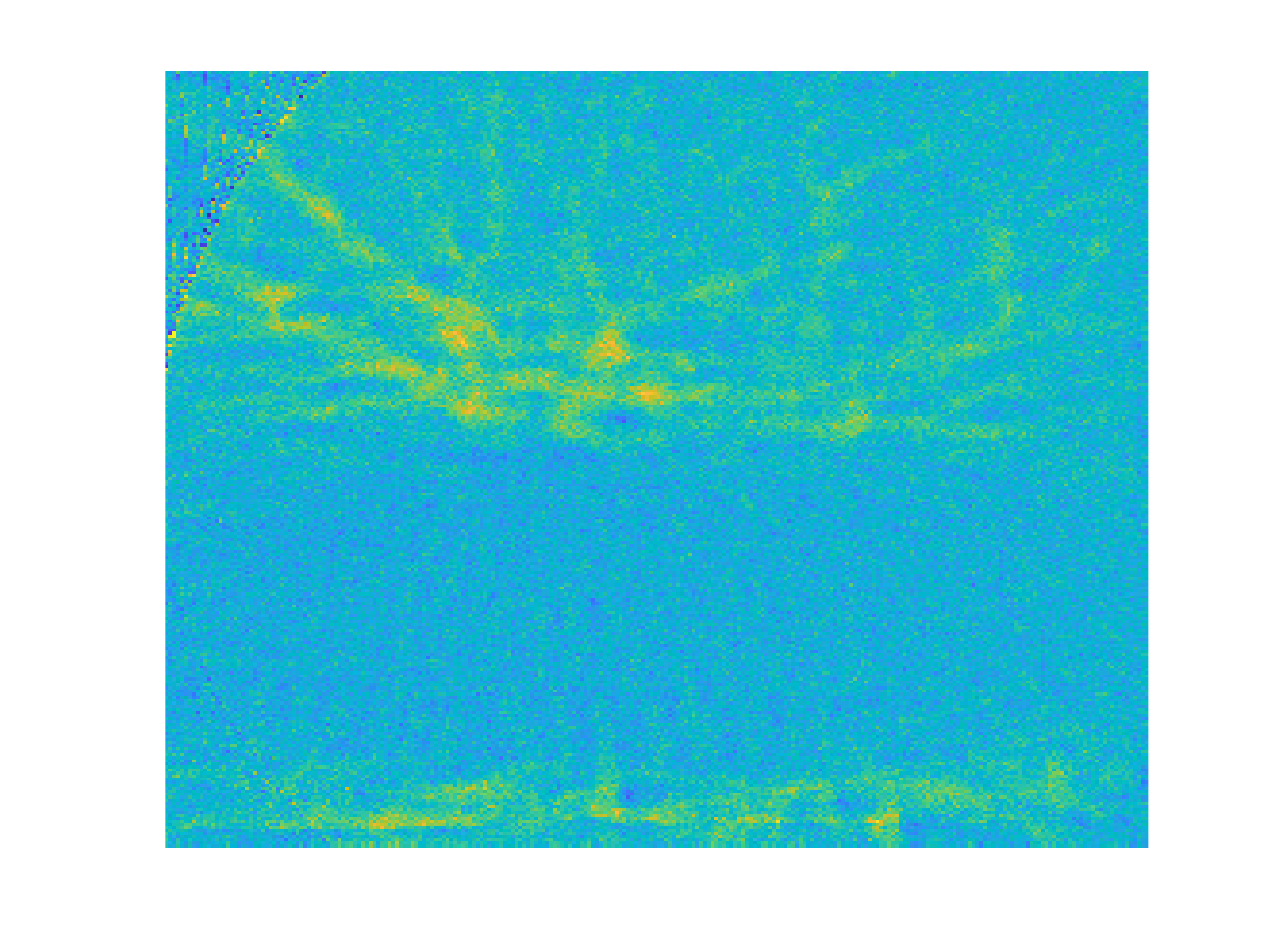}
    \includegraphics[width=0.3\columnwidth, height=0.25\columnwidth]{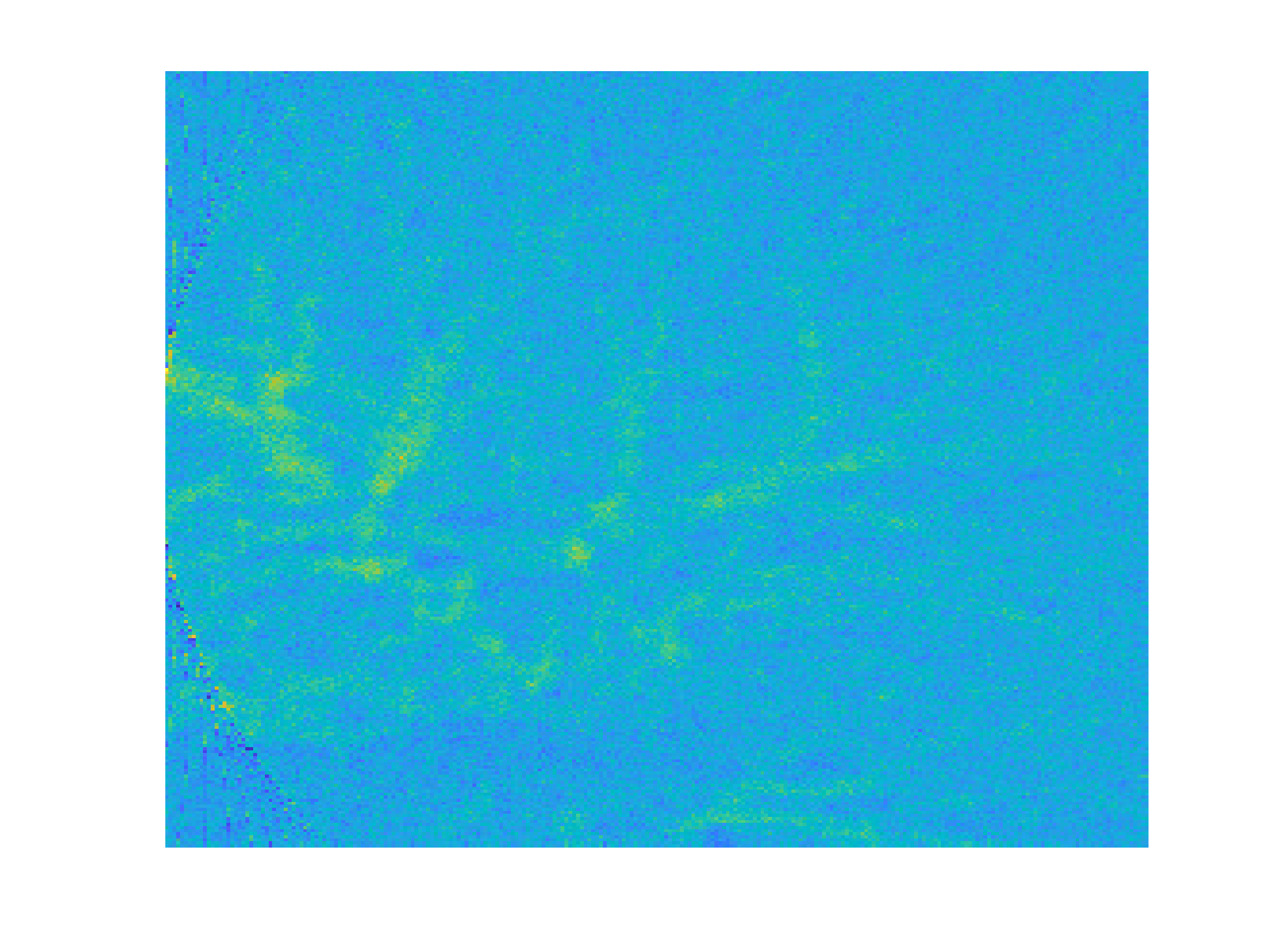}\\
    \includegraphics[width=0.3\columnwidth, height=0.25\columnwidth]{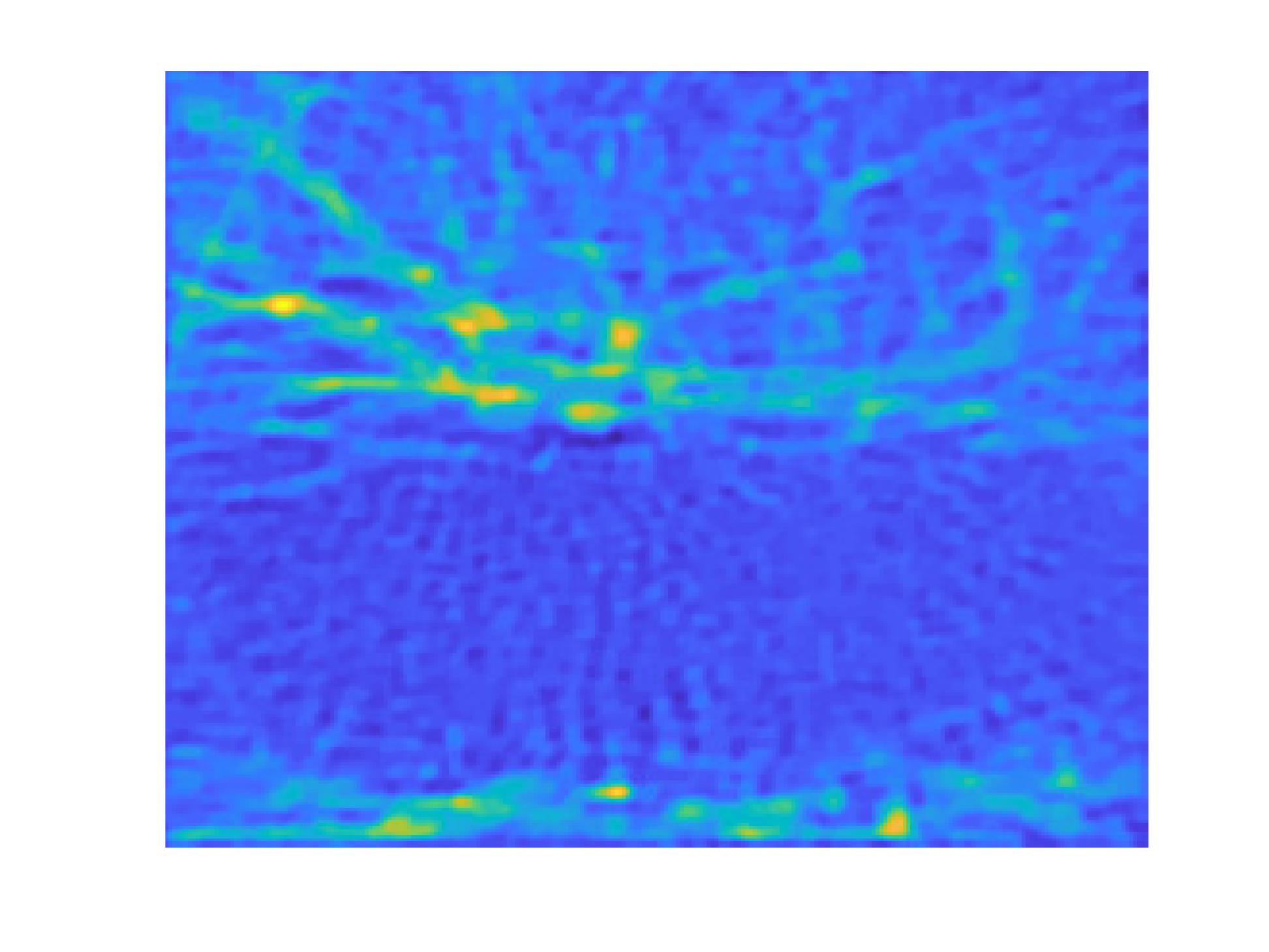}
    \includegraphics[width=0.3\columnwidth, height=0.25\columnwidth]{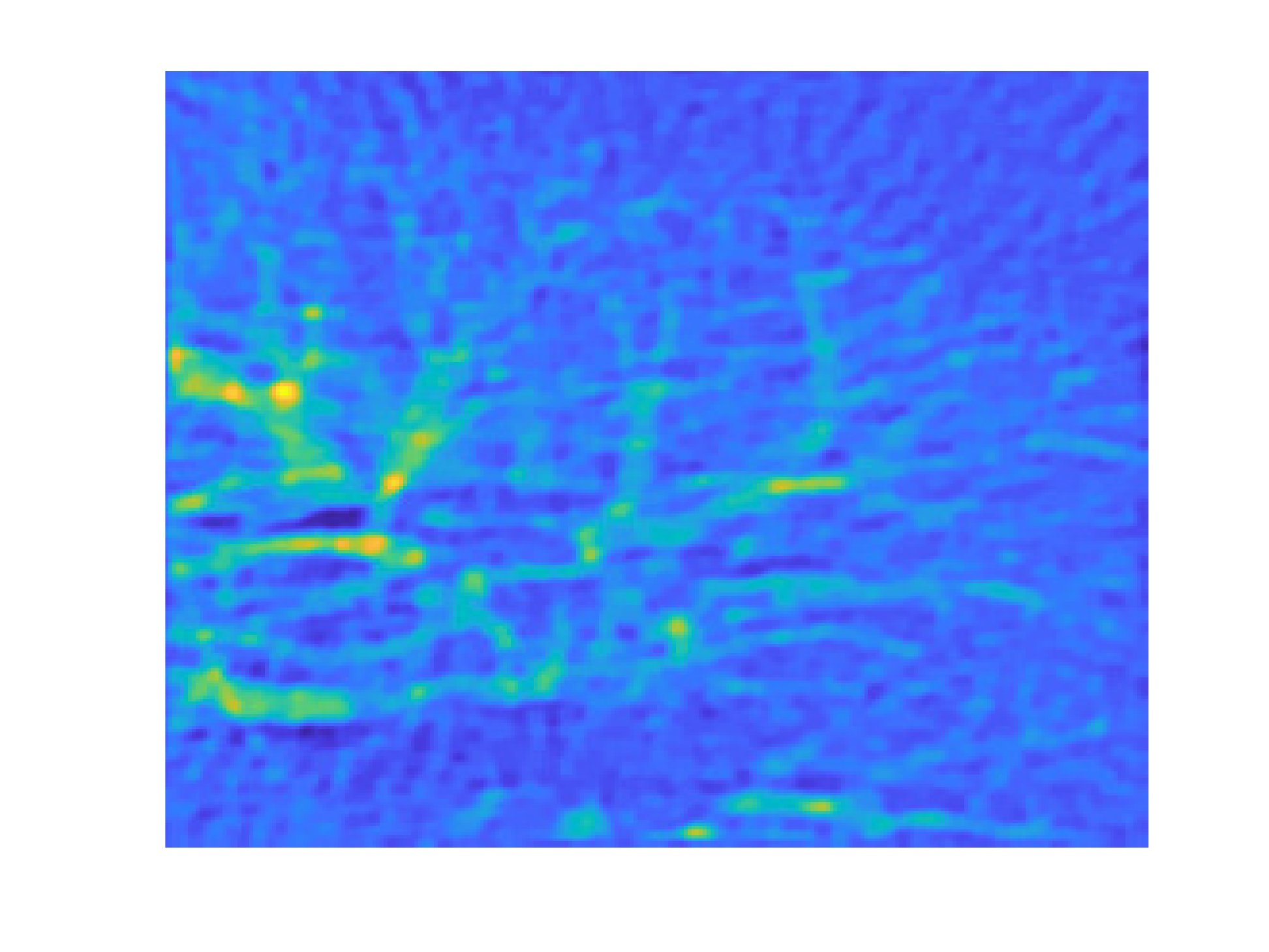}\\
    \includegraphics[width=0.3\columnwidth, height=0.25\columnwidth]{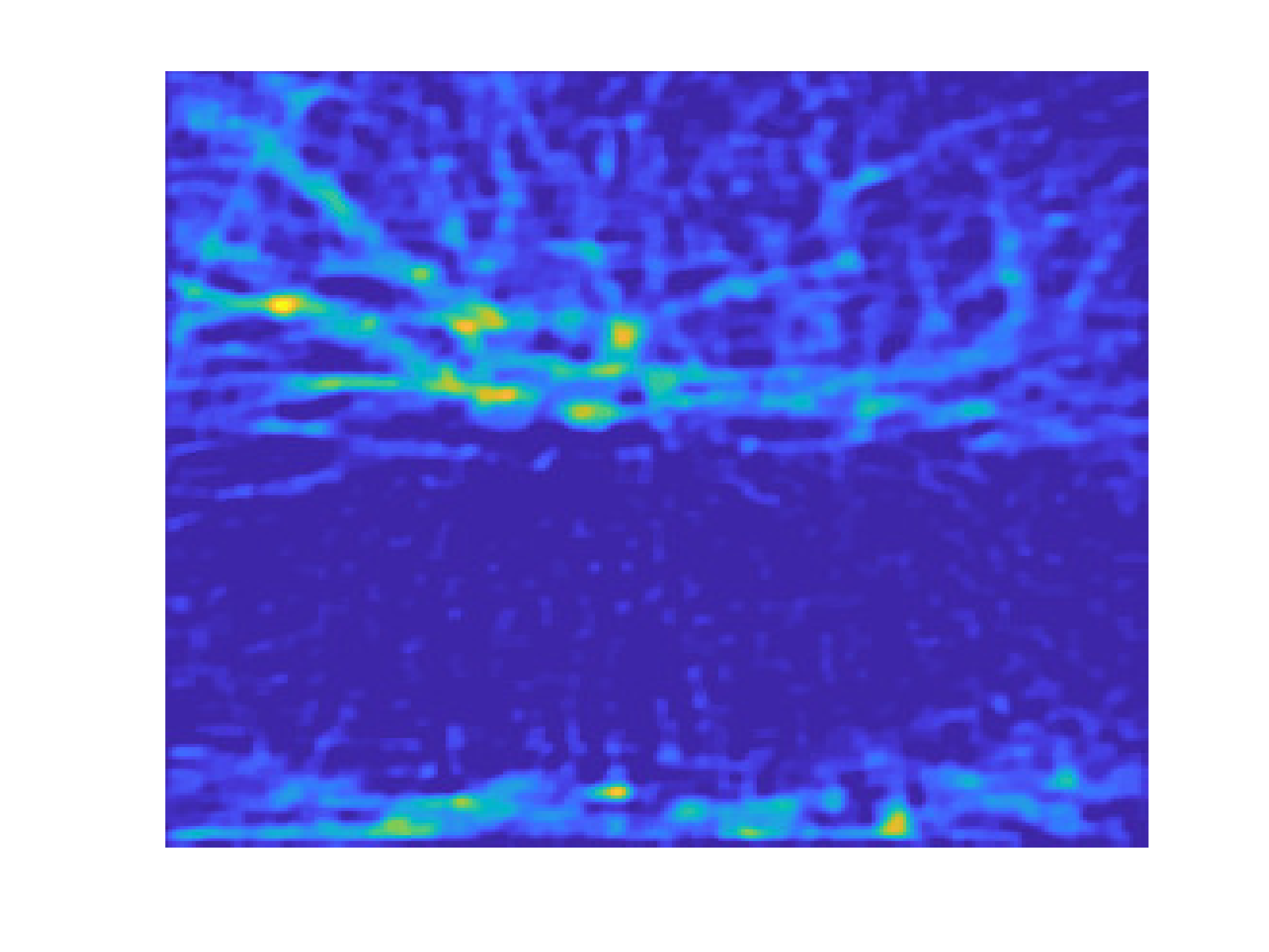}
    \includegraphics[width=0.3\columnwidth, height=0.25\columnwidth]{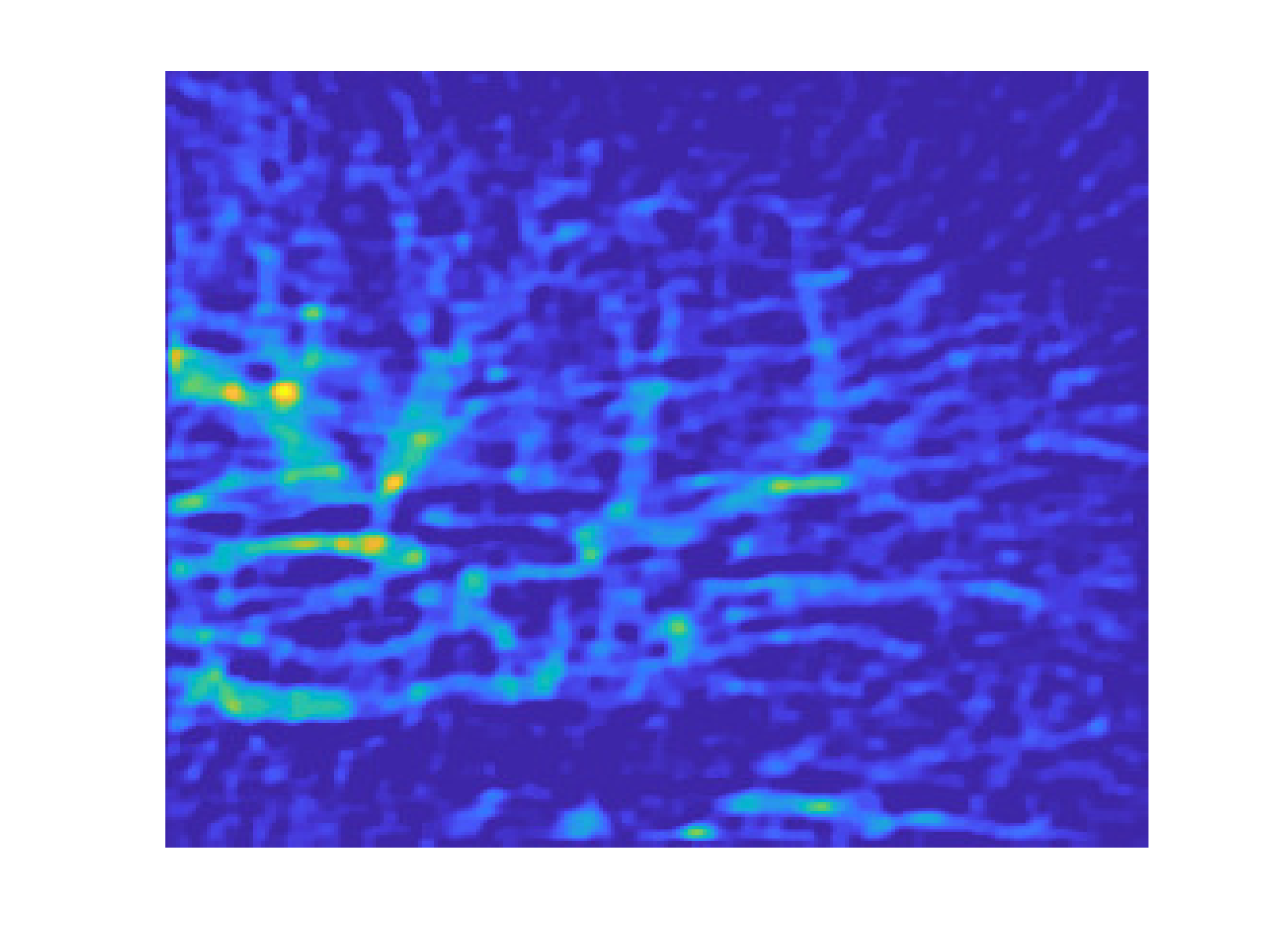}\\
    \includegraphics[width=0.3\columnwidth, height=0.25\columnwidth]{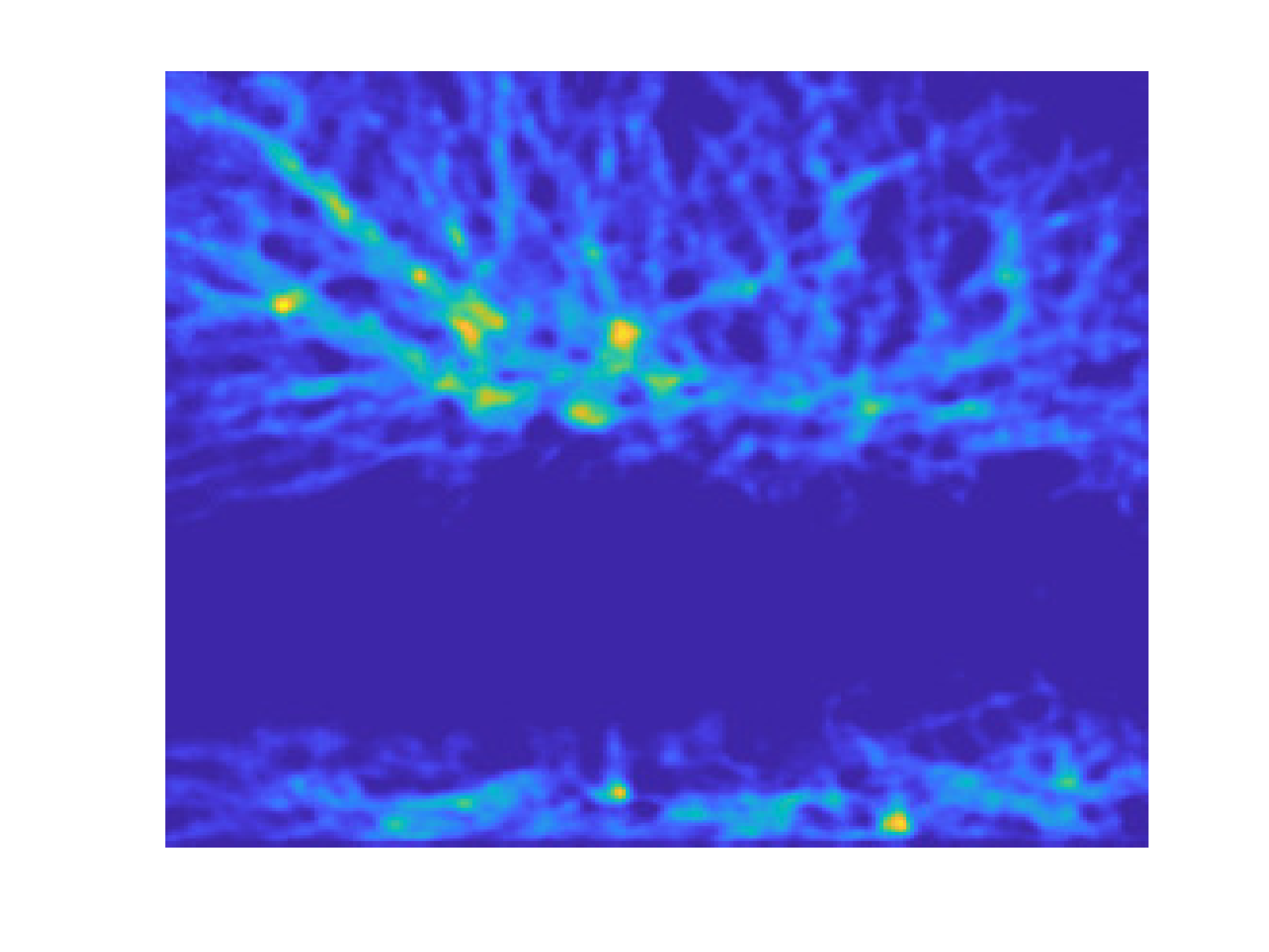}
    \includegraphics[width=0.3\columnwidth, height=0.25\columnwidth]{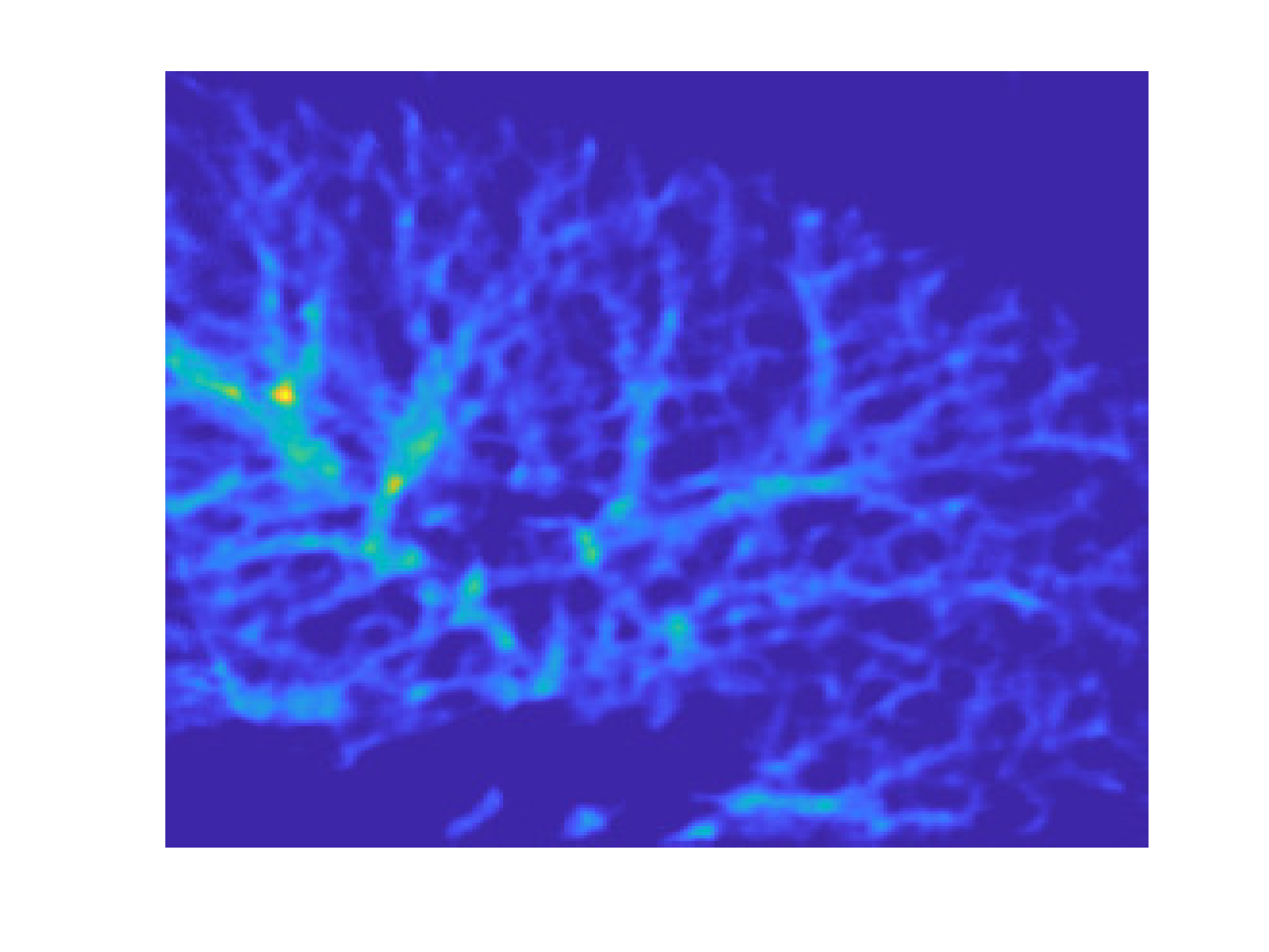}
        \caption{\textbf{Reconstructions from (noisy) evaluation data:}\label{fig:noisy}
    First Row: Projection images of lung vessel phantom.
    {Second Row:} DAL reconstructions.
    {Third Row:}     TV-minimization.
    {Fourth Row:}   TV-minimization with positivity constraint.
    {Bottom Row:}  Proposed DALnet.
     The reconstruction region is $[\SI{-12.5}{mm}, \SI{12.5}{mm}] \times [\SI{-20}{mm}, \SI{5}{mm}]$ and the detection curve is $\rand = \set{\rrs \colon \norm{\rrs}_2 = \SI{50}{mm} \wedge \rrs_2 < 0}$.}
\end{figure}

\subsection{Simulation results}

Figure~\ref{fig:noisy} shows a comparison of the reconstruction  algorithms applied to 
two randomly chosen examples from evaluation data.  The  computation time
recover a single PA projection image with TV minimization using  30 iterations  
implemented in Matlab is $\SI{48}{seconds}$ on a AMD Ryzen 7 1700X CPU.
The DALnet is implemented in Keras \cite{chollet2015keras}
 on top of Tensorflow  \url{https://www.tensorflow.org}.
The training was done on a NVIDIA TITAN Xp GPU and lasted for about
$\SI{10}{hours}$ for the whole training procedure. The application of the DALnet 
(including the learned UBP and the Unet) to one data set only requires 
$\SI{0.018}{seconds}$, and evaluation of the UBP requires $\SI{0.013}{seconds}$.
This yields  to a frame rate of more than 50 reconstructed PA images per second,  
clearly allowing PA real-time monitoring.    

All reconstructions with the UBP show  typical angular stripe-like under-sampling artifacts;  additionally the results are blurred according to the system PSF. 
The results for TV-minimization  with positivity constraint   and  the proposed DALnet yield
reconstructions almost  free from under-sampling artifacts and outperform the other reconstruction  methods.     Additionally, they are  capable to remove (at least partially) the  blurring due the the PSF.
DALnet even outperforms positivity constraint TV in terms of visual image quality.  
For TV regularization the regularization  parameter has been adjusted manually to yield visually appealing results. The ill-posedness of the involved convolution prohibits decreasing the regularization parameter as otherwise the noise turned out to be severely amplified in the reconstruction. 
The quite large  
regularization parameter, however,  yields  to over-smoothing of the fine
blood-vessel structures with TV-minimization. The DALnet does not suffer from
this limitation and  yields high-resolution images without noise amplification.

\subsection{Error analysis}

For a more quantitative evaluation of the performance of the above reconstruction
methods,  we evaluate various error measures  averaged  over  the  200 test images
not contained in the training data.  We evaluate the reconstructions in terms of 
scaled and shifted relative $\ell^2$ and  $\ell^1$ distance,  the scaled and shifted  structured similarity index
(SSIM; a common  measure  for predicting the perceived quality of digital images introduced in \cite{wang2004image}),  and the correlation.  
Following   \cite{jin2017deep,hauptmann2018model}, the
scaled and shifted relative $\ell^2$ and  $\ell^1$ distance and  scaled and shifted  
 SSIM between two images $\HH, \XX  \in \R^{\Nx \times \Nx}$ are computed by   
\begin{align*}
\mathcal{E}_{1}(  \HH, \XX ) 
& \coloneqq \min_{\alpha,\beta \in \R} \frac{\norm{ \alpha  \HH - \XX - \beta}_1}{\norm{ \XX }_1} \,, \\
\mathcal{E}_{2}(  \HH, \XX ) 
& \coloneqq \min_{\alpha,\beta \in \R} \frac{\norm{ \alpha  \HH - \XX - \beta}_2}{\norm{ \XX }_2} \,, \\
\mathcal{E}_{\operatorname{SSIM}} (\HH, \XX)
& \coloneqq \max_{\alpha,\beta\in \R}\{ \operatorname{SSIM}(\alpha \HH - \beta, \XX)\} \,,
\end{align*}
respectively.  Here $\enorm{}_{p}$ is the standard  $\ell^p$ distance and $\operatorname{SSIM}$ is   the structured similarity index. Note that a larger value of $\mathcal{E}_{\operatorname{SSIM}} (\HH, \XX)$  corresponds to a higher structural similarity of $\HH$ and $\XX$.  
Opposed to that,  smaller values of $\mathcal{E}_{2}$, $\mathcal{E}_{1}$ correspond to  
more similar results.  

Table~\ref{tab:err} shows the computed  error measures  for reconstructions with plain UBP, TV, TV with positivity constraint and the proposed DALnet, averaged over  200 test images
not contained in the training set. It can be observed that in terms of all error measures DALnet outperforms the other reconstruction measures. Also notably, including positivity significantly 
improves  TV reconstruction in terms of all the quality measures.

\begin{table}
\centering
\begin{tabular}{ l | c  c  c  c}
\toprule
error measure  & UBP & TV & TV pos  & DALnet  \\
\midrule
    $\ell^2$       & 0.555 & 0.112 & 0.102  & 0.085 \\ 
    $\ell^1$       & 0.860 & 0.418 & 0.361  & 0.319 \\
    SSIM        & 0.305 & 0.586 & 0.678  & 0.726 \\ 
    correlation & 0.382 & 0.911 & 0.919  & 0.933 \\
\bottomrule
\end{tabular}
\caption{\textbf{Various error measures averaged over 200 test images not 
contained in the training set.}\label{tab:err}
We evaluate the reconstructions in terms of  scaled and shifted relative  $\ell^2$ and  $\ell^1$  distances,  the scaled and shifted  SSIM,  and the correlation.}
\end{table}

\begin{figure}[htb!]
    \centering
    \includegraphics[width=0.3\columnwidth, height=0.25\columnwidth]{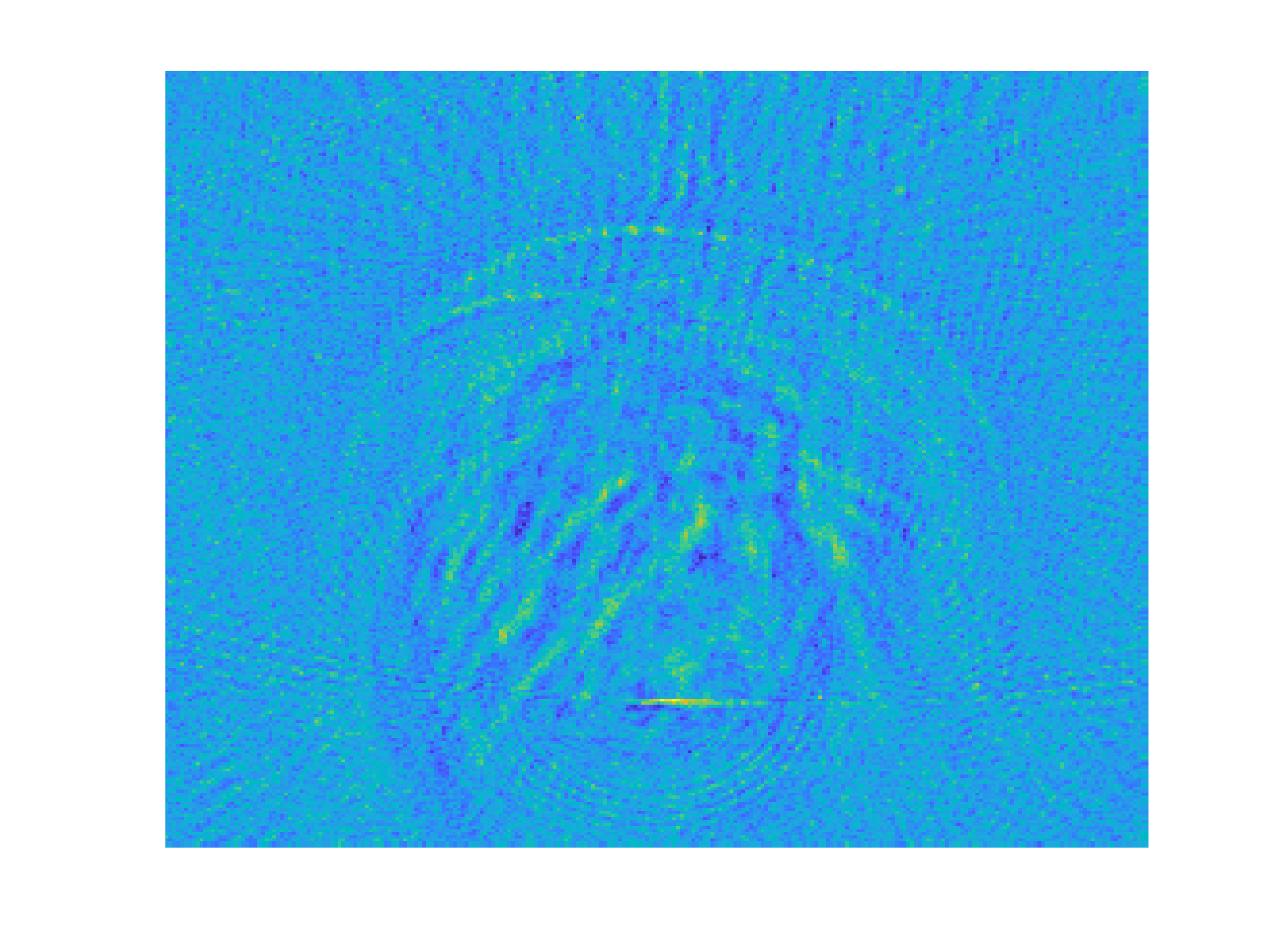}
    \includegraphics[width=0.3\columnwidth, height=0.25\columnwidth]{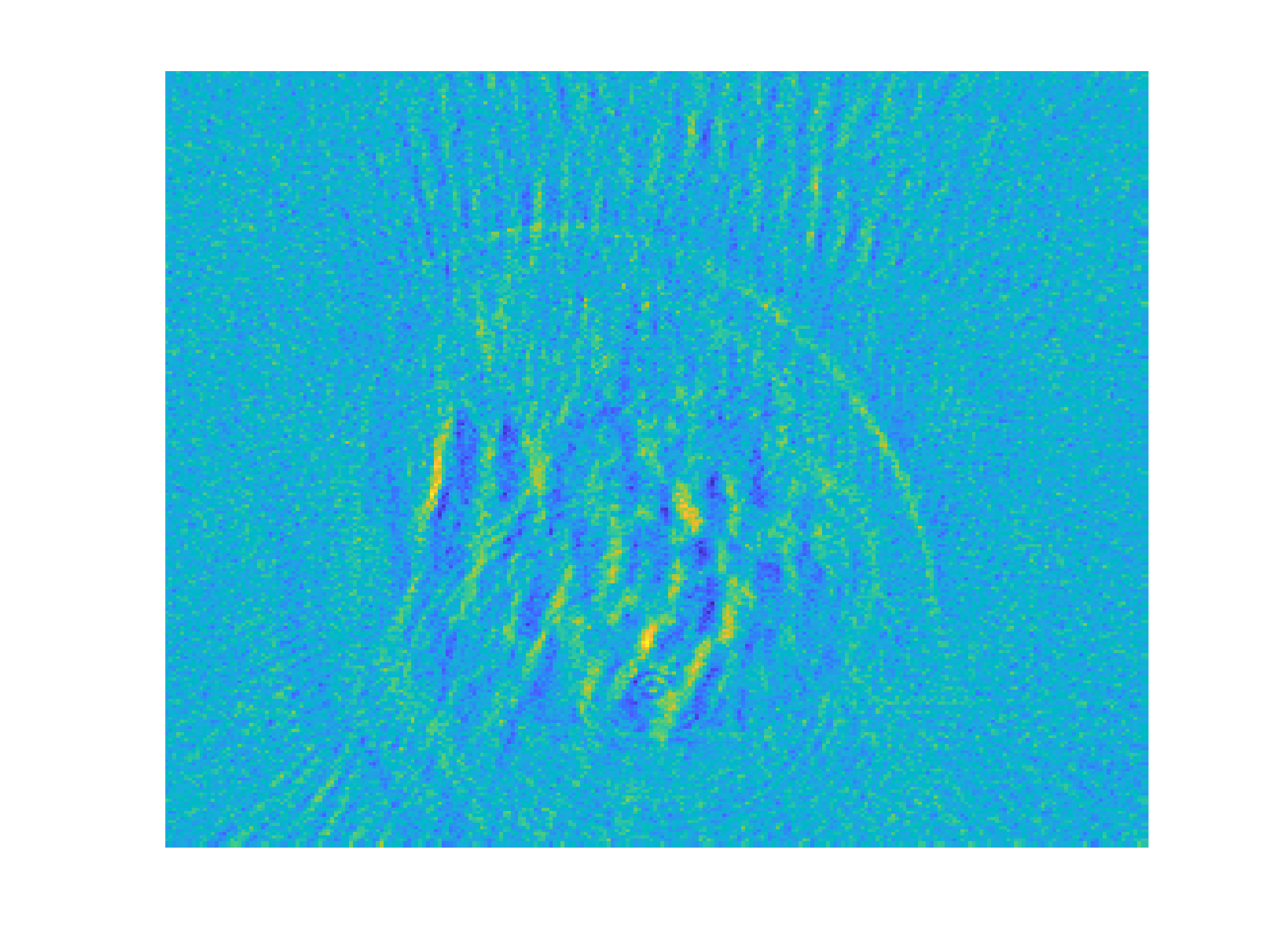}\\
    \includegraphics[width=0.3\columnwidth, height=0.25\columnwidth]{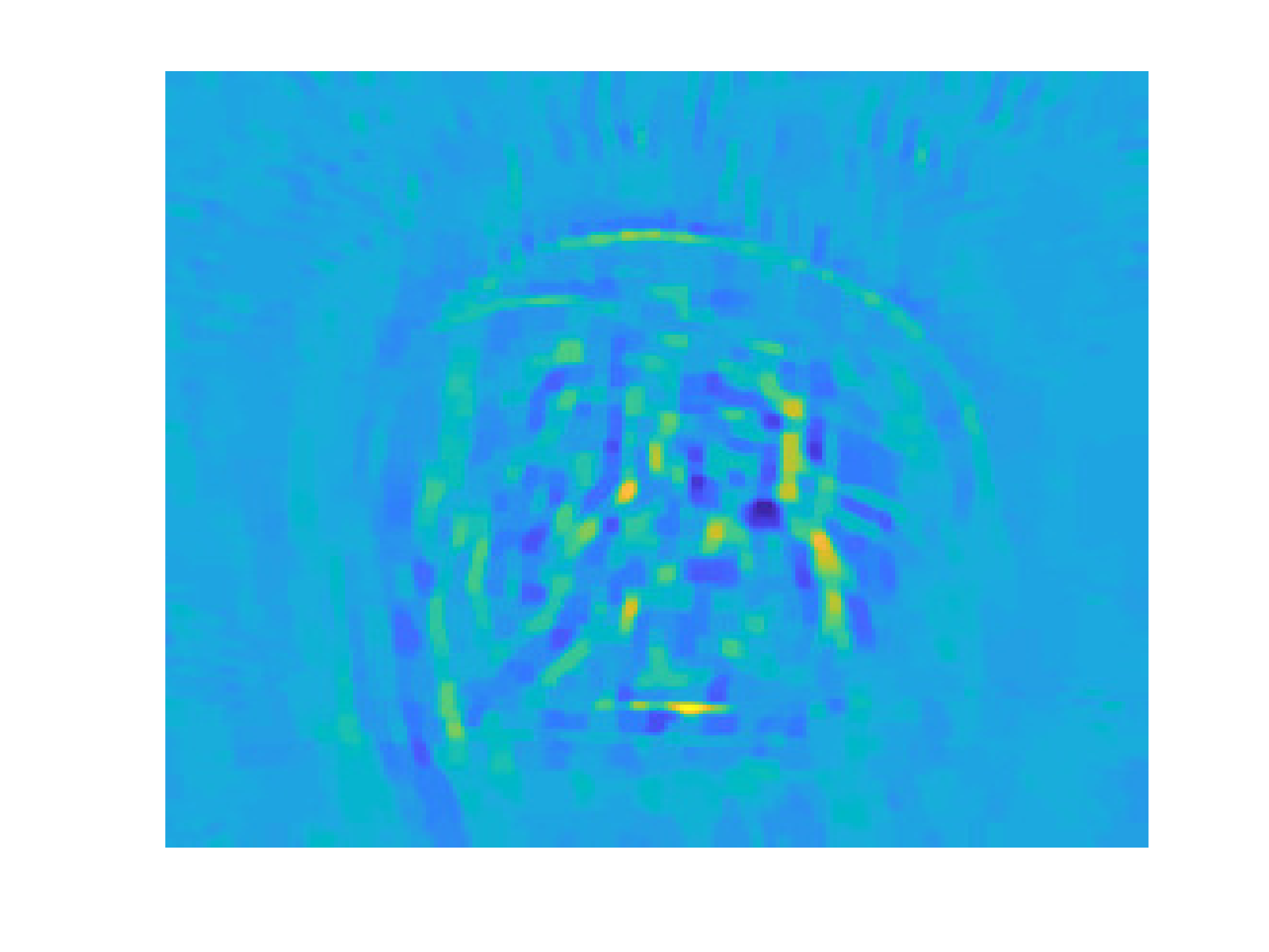}
    \includegraphics[width=0.3\columnwidth, height=0.25\columnwidth]{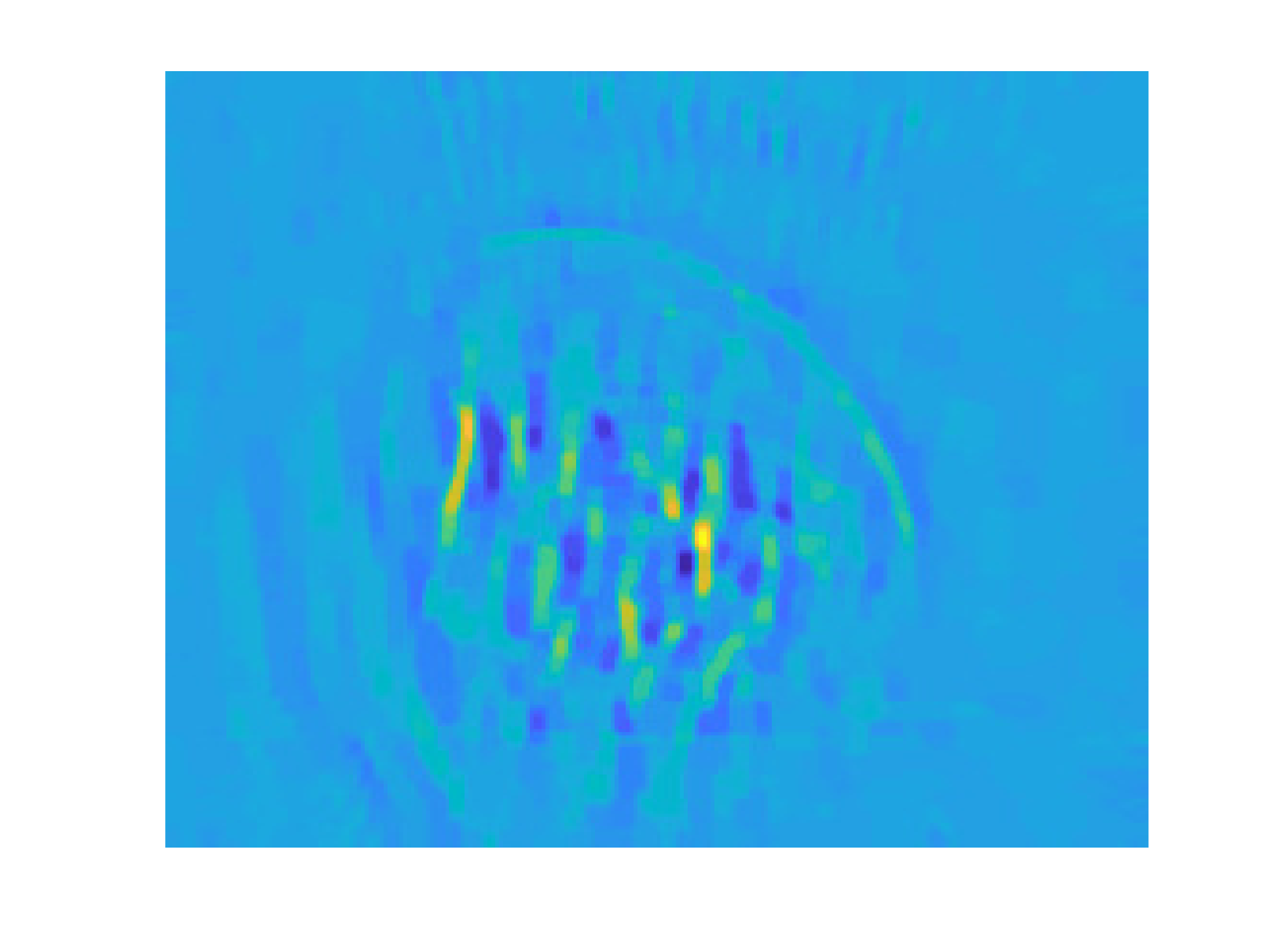}\\
    \includegraphics[width=0.3\columnwidth, height=0.25\columnwidth]{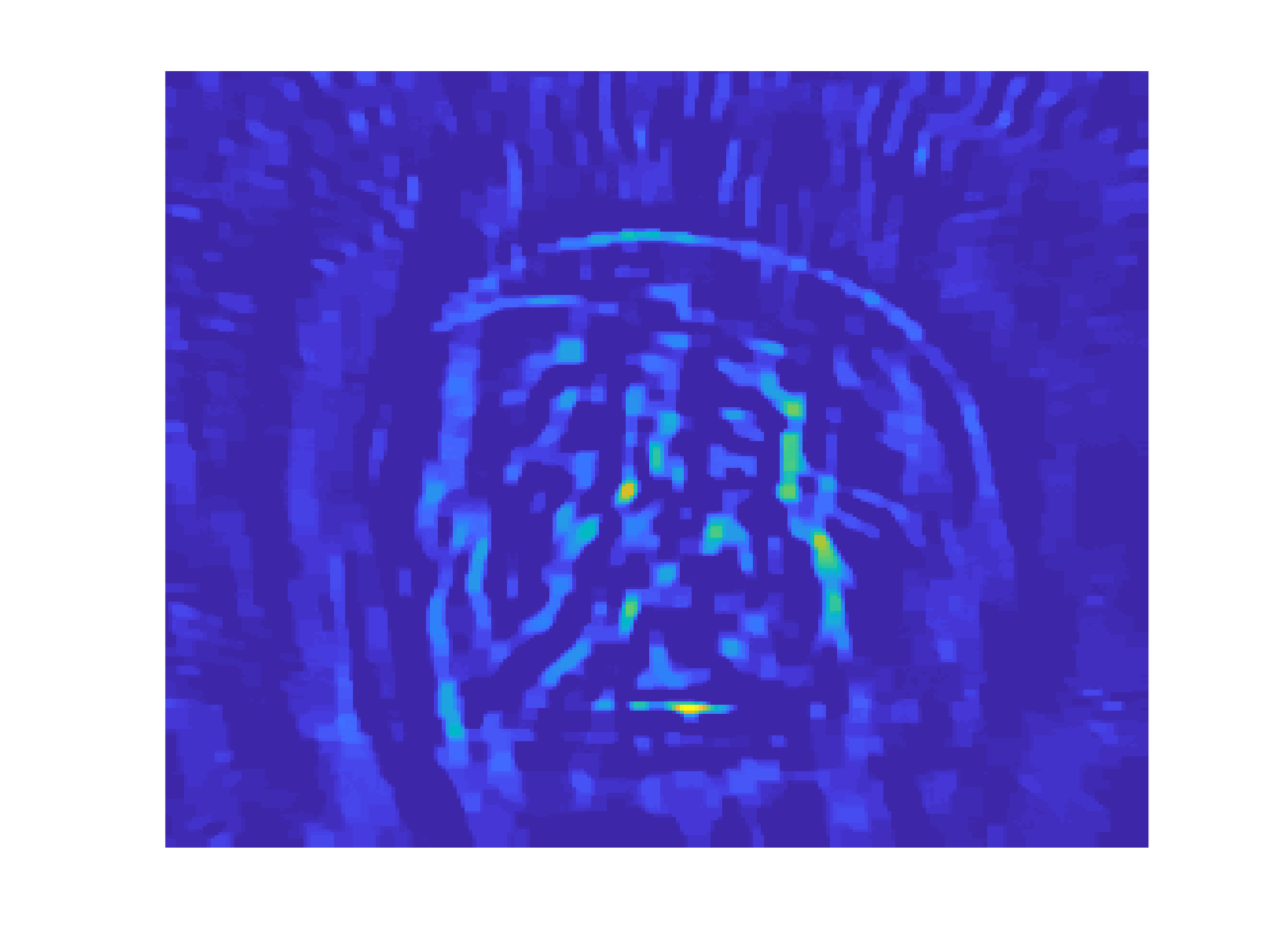}
    \includegraphics[width=0.3\columnwidth, height=0.25\columnwidth]{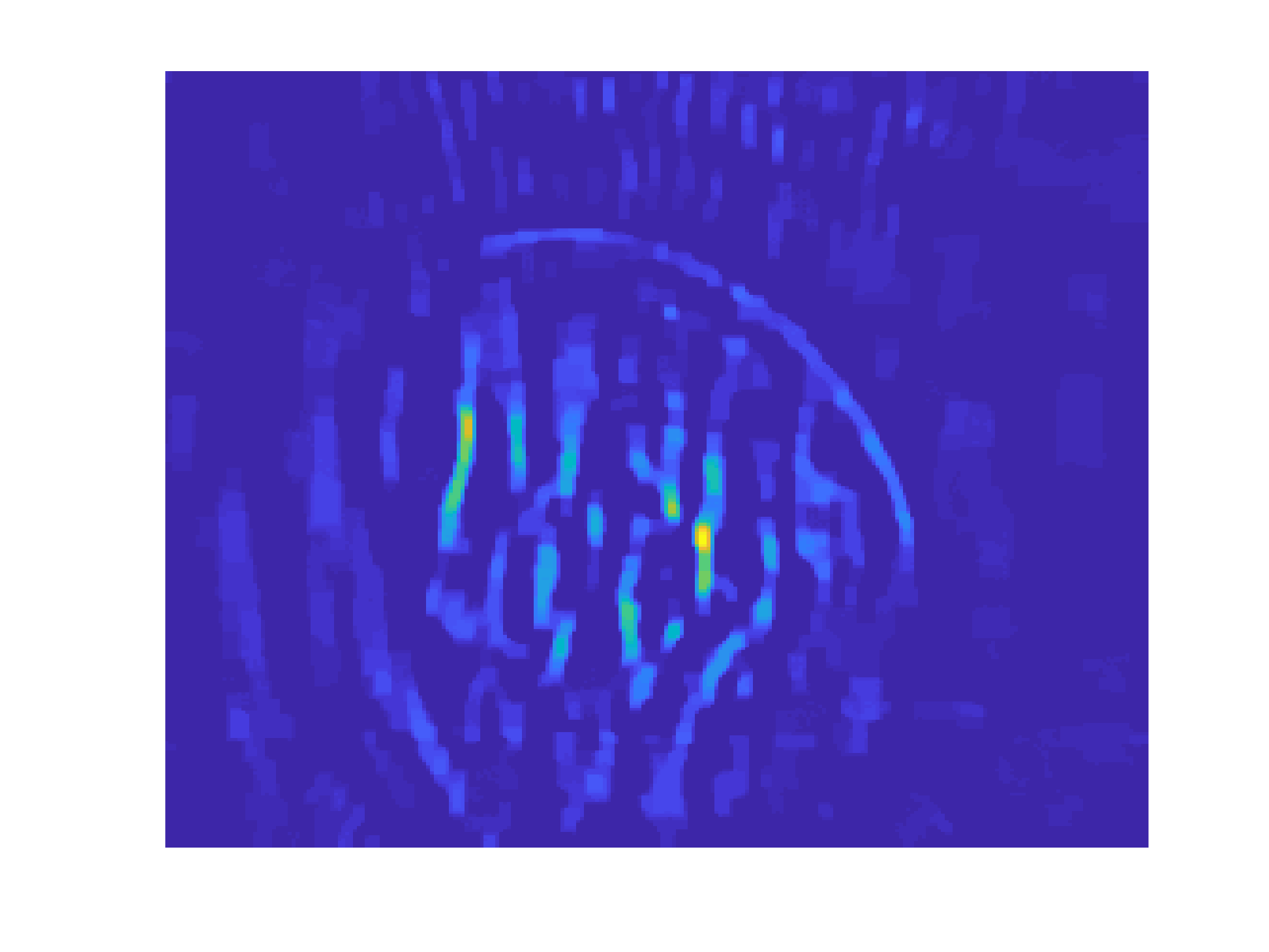}\\
    \includegraphics[width=0.3\columnwidth, height=0.25\columnwidth]{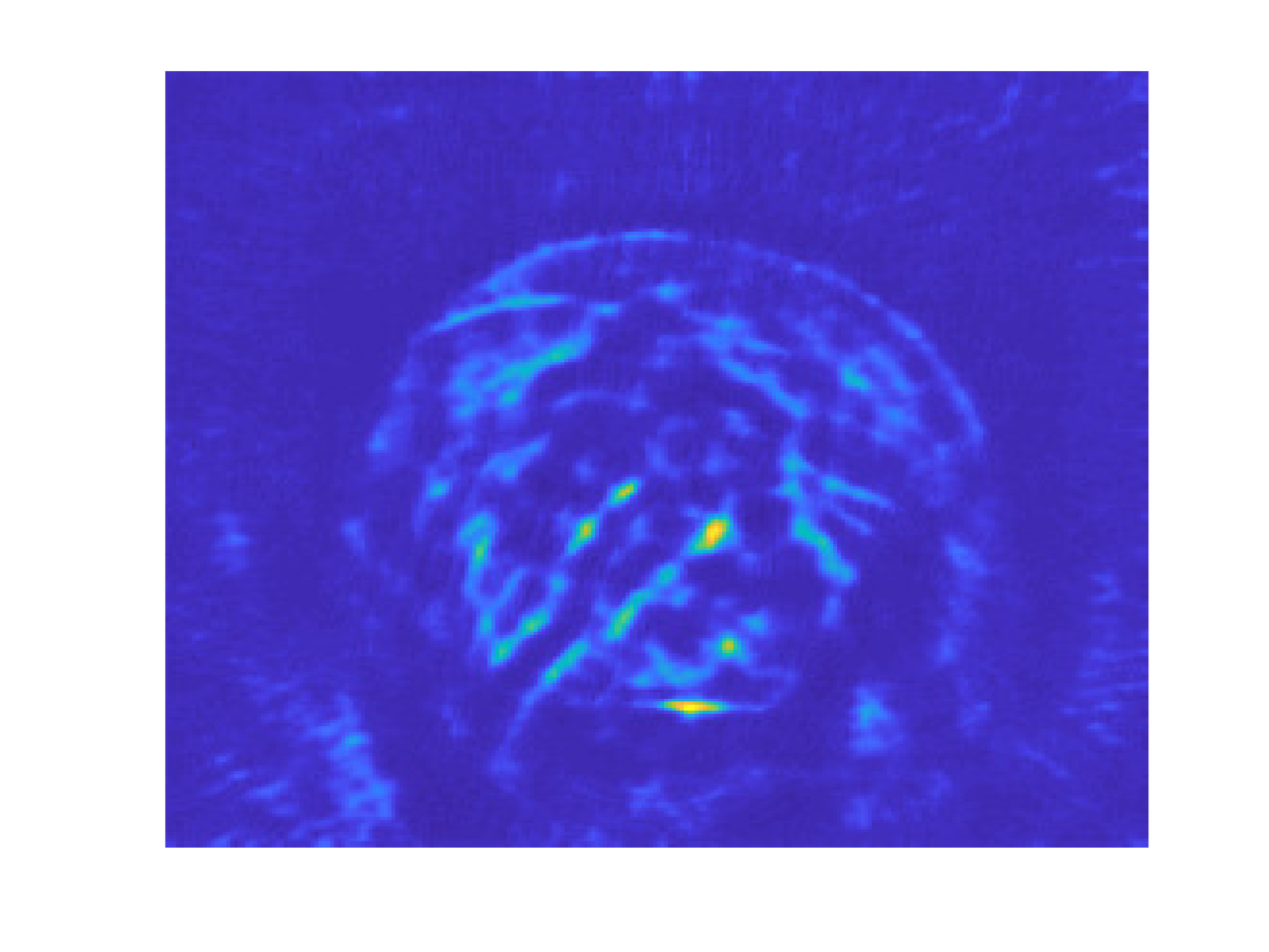}
    \includegraphics[width=0.3\columnwidth, height=0.25\columnwidth]{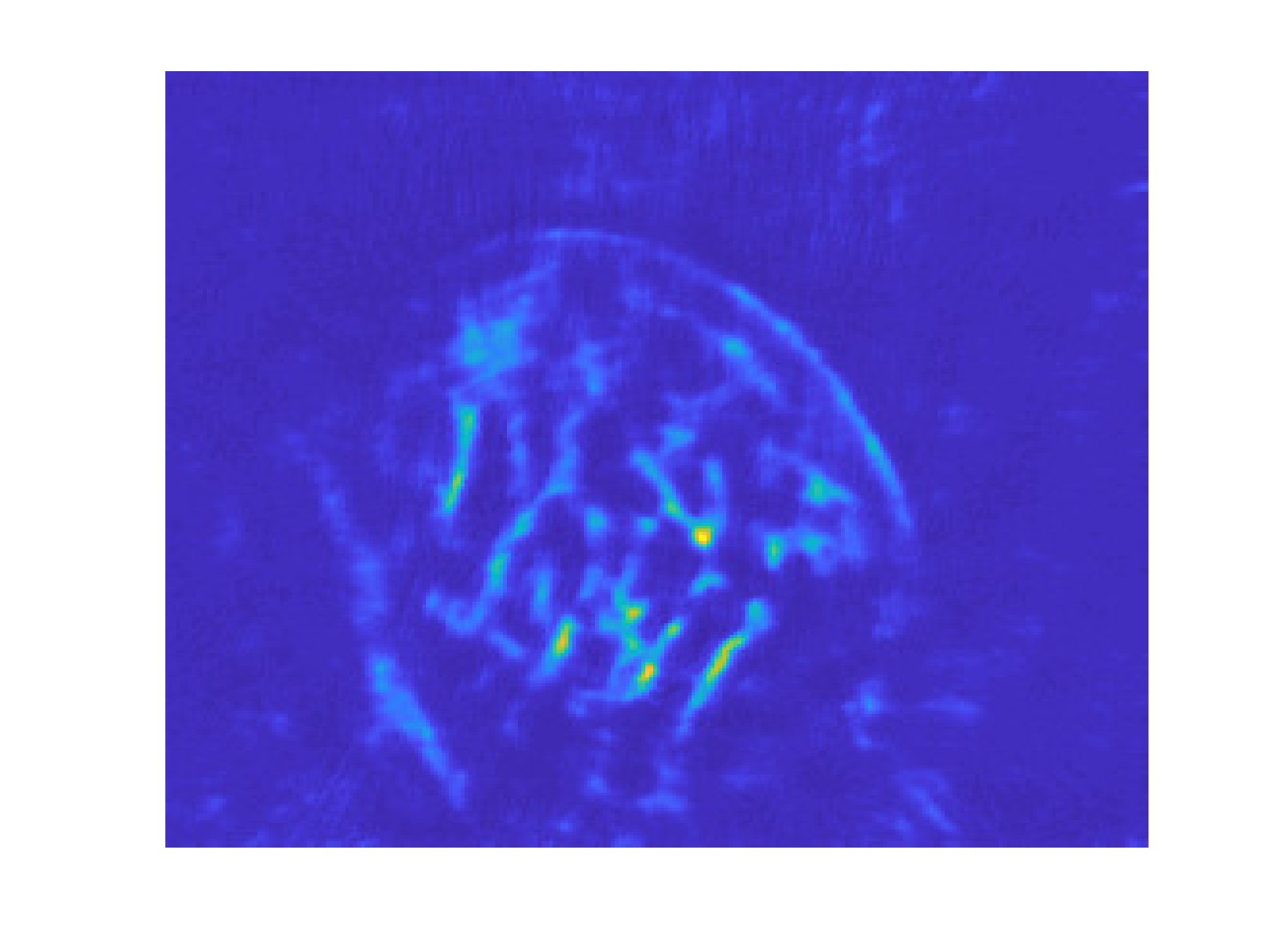}
        \caption{\textbf{Reconstructions from  experimental data:} \label{fig:finger}
        {Top Row:} UBP reconstructions.
       {Second Row:}     TV-minimization.
      {Third  Row:}   TV-minimization with positivity constraint.
       {Bottom Row:}  Reconstructions using proposed DALnet.
        The reconstructions correspond to three different rotation angles; 
        the reconstruction region is 
        $[\SI{-12.5}{mm}, \SI{12.5}{mm}] \times  [\SI{-20}{mm}, \SI{5}{mm}]$; the detectors curve
        is $\rand = \set{\rrs \colon \norm{\rrs}_2 = \SI{50}{mm} \wedge \rrs_2 < 0}$.}
        \end{figure}

\subsection{Results for in-vivo data}

We also evaluate  the reconstruction algorithms on experimental data collected from the 64-line PA system.
We therefore have taken  several snapshots of a finger of one of the authors from  200  different rotation angles
around a single axis, and apply the  DALnet trained above, the UBP  and TV minimization.  For TV minimization, we use 30 iterations of the Chambolle-Pock algorithm with  the same settings as in the noisy simulated data case. For noise reduction, the PA signals of the human finger are four times averaged prior to the reconstruction process.
Two randomly selected   reconstructed PA projection images  are  shown in  Figure~\ref{fig:finger}. The  reconstructions are again evaluated at  $256 \times 256$ regular sampling points in the square $[\SI{-12.5}{mm}, \SI{12.5}{mm}] \times  [\SI{-20}{mm}, \SI{5}{mm}]$ with the  detector curve $\rand = \set{\rrs \colon \norm{\rrs}_2 = \SI{50}{mm} \wedge \rrs_2 < 0}$.

From Figure~\ref{fig:finger} one observes that the UBP  reconstruction  of the finger suffers
from streak like under-sampling-artifacts as well as blurring due to the PSF of the PA imaging
system. Both,  TV minimization and the DALnet are able to partially remove
these artifacts.  The reconstructions using the DALnet are sharper than the
TV reconstructions (with or without positivity constraint) and
better capable of removing blurring due to the PSF.
In the TV reconstruction, any selection of the regularization parameter $\la$ is a 
tradeoff between amplifying error and resolving the fine structures of the blood vessels inside the
finger. The reconstructions using the DALnet do not suffer from such a trade-off.
For  TV-minimization  (and related variational or iterative algorithms),
blood vessel structures and  under-sampling artifacts are  hardly distinguishable.
Increasing the regularization parameter and thereby putting
more emphasis on artifact (and noise) removal at the same time also removes fine structures
in the blood vessels. Opposed to that, the  DALnet is trained  to distinguish between blood vessels
and under-sampling artifacts and therefore is capable of producing high-resolution reconstructions  while also removing under-sampling artifacts.

Following the standard procedure for PAT with integrating line detectors  \cite{BurBauGruHalPal07,haltmeier2009frequency,paltauf2009photoacoustic}, 
the PA  images  have been converted to a 3D image $\mathbf{f}_{\rm rec} \in \R^{\Nx \times \Nx \times \Nx}$ of a finger. Supplementary videos 1-4 (see \url{https://applied-math.uibk.ac.at/cms/index.php/preprints-2018}) 
visualize the 3D reconstructions for the FBP, TV with and without positivity constraint and the proposed DALnet. 
From these videos DALnet seems to produce the best three dimensional reconstructions of Roberts finger.

\section{Conclusion}
\label{sec:discussion}

In this study we reported a deep learning approach for high-resolution
PA projection imaging. The proposed  DALnet seems to be the
first network for PAT that is capable to account for limited view,
under-sampling artifacts and blurring due the system PSF.
We designed and trained the DALnet for a piezoelectric 64-line detector system
\cite{paltauf2017piezoelectric}. The obtained direct reconstruction
algorithm requires less that  1/50  seconds   on a  standard desktop  PC
in with  an NVIDIA TITAN Xp GPU  (available for about  \SI{1400}{euro}) 
to reconstruct  a $256 \times 256$ PA projection image.   We are not aware of  
any PAT system that is capable to produce such high-resolution real-time 
monitoring of 3D structures at such low hardware costs.

Future work will be done to compare the proposed DALnet with other deep learning image reconstruction approaches including the null-space network~\cite{schwab2018deep}, iterative networks~\cite{adler2017solving,kelly2017deep}, learned  projected gradient~\cite{gupta2018cnn}, variational networks~\cite{kobler2017variational} and  the learned network regularizers~\cite{li2018nett}. Particular emphasis will be given on the evaluation on real world   data.
Finally, we point out that the  proposed network is not restricted to PA imaging.
The design of  similar networks for other  tomographic modalities is possible and  
subject of future investigations. This includes fan beam  and spiral CT, where weight 
corrections  are inherently included in FPB algorithms that are still standard in commercial CT 
scanners.

\section*{Acknowledgement}

SA and MH  acknowledge support of the Austrian Science Fund (FWF), project P 30747-N32. The work of RN has been supported by the FWF, project P 28032. JS, SA and MH  acknowledge
 support of NVIDIA Corporation with the donation of the Titan Xp GPU used for this research.


\begin{thebibliography}{10}
\providecommand{\url}[1]{#1}
\csname url@samestyle\endcsname
\providecommand{\newblock}{\relax}
\providecommand{\bibinfo}[2]{#2}
\providecommand{\BIBentrySTDinterwordspacing}{\spaceskip=0pt\relax}
\providecommand{\BIBentryALTinterwordstretchfactor}{4}
\providecommand{\BIBentryALTinterwordspacing}{\spaceskip=\fontdimen2\font plus
\BIBentryALTinterwordstretchfactor\fontdimen3\font minus
  \fontdimen4\font\relax}
\providecommand{\BIBforeignlanguage}[2]{{%
\expandafter\ifx\csname l@#1\endcsname\relax
\typeout{** WARNING: IEEEtran.bst: No hyphenation pattern has been}%
\typeout{** loaded for the language `#1'. Using the pattern for}%
\typeout{** the default language instead.}%
\else
\language=\csname l@#1\endcsname
\fi
#2}}
\providecommand{\BIBdecl}{\relax}
\BIBdecl

\bibitem{Bea11}
P.~Beard, ``Biomedical photoacoustic imaging,'' \emph{Interface focus}, vol.~1,
  no.~4, pp. 602--631, 2011.

\bibitem{Wan09b}
L.~V. Wang, ``Multiscale photoacoustic microscopy and computed tomography,''
  \emph{Nature Phot.}, vol.~3, no.~9, pp. 503--509, 2009.

\bibitem{BurHofPalHalSch05}
P.~Burgholzer, C.~Hofer, G.~Paltauf, M.~Haltmeier, and O.~Scherzer,
  ``Thermo\-acoustic tomography with integrating area and line detectors,''
  \emph{IEEE Trans. Ultrasonic and Frequency Control}, vol.~52, no.~9, pp.
  1577--1583, September 2005.

\bibitem{PalNusHalBur07a}
G.~Paltauf, R.~Nuster, M.~Haltmeier, and P.~Burgholzer, ``Photoacoustic
  tomography using a {M}ach-{Z}ehnder interferometer as an acoustic line
  detector,'' \emph{App. Opt.}, vol.~46, no.~16, pp. 3352--3358, 2007.

\bibitem{BurBauGruHalPal07}
P.~Burgholzer, J.~Bauer-Marschallinger, H.~Gr{\"u}n, M.~Haltmeier, and
  G.~Paltauf, ``Temporal back-projection algorithms for photoacoustic
  tomography with integrating line detectors,'' \emph{Inverse Probl.}, vol.~23,
  no.~6, pp. S65--S80, 2007.

\bibitem{haltmeier2009frequency}
M.~Haltmeier, ``Frequency domain reconstruction for photo- and thermo\-acoustic
  tomography with line detectors,'' \emph{Math. Models Methods Appl. Sci.},
  vol.~19, no.~2, pp. 283--306, 2009.

\bibitem{paltauf2009photoacoustic}
G.~Paltauf, R.~Nuster, M.~Haltmeier, and P.~Burgholzer, ``Photoacoustic
  tomography with integrating area and line detectors,'' in \emph{Photoacoustic
  imaging and spectroscopy}.\hskip 1em plus 0.5em minus 0.4em\relax CRC Press,
  2009, ch.~20, pp. 251--263.

\bibitem{GraEtAl14}
S.~Gratt, R.~Nuster, G.~Wurzinger, M.~Bugl, and G.~Paltauf, ``64-line-sensor
  array: fast imaging system for photoacoustic tomography,'' \emph{Proc. SPIE},
  vol. 8943, p. 894365, 2014.

\bibitem{BauPSP15}
J.~Bauer-Marschallinger, K.~Felbermayer, K.-D. Bouchal, I.~A. Veres,
  H.~Gr{\"u}n, P.~Burgholzer, and T.~Berer, ``Photoacoustic projection imaging
  using a 64-channel fiber optic detector array,'' in \emph{Proc. SPIE}, vol.
  9323, 2015.

\bibitem{paltauf2017piezoelectric}
G.~Paltauf, P.~Hartmair, G.~Kovachev, and R.~Nuster, ``Piezoelectric line
  detector array for photoacoustic tomography,'' \emph{Photoacoustics}, vol.~8,
  pp. 28--36, 2017.

\bibitem{Bauer-Marschallinger:17}
J.~Bauer-Marschallinger, K.~Felbermayer, and T.~Berer, ``All-optical
  photoacoustic projection imaging,'' \emph{Biomed. Opt. Express}, vol.~8,
  no.~9, pp. 3938--3951, 2017.

\bibitem{arridge2016accelerated}
S.~Arridge, P.~Beard, M.~Betcke, B.~Cox, N.~Huynh, F.~Lucka, O.~Ogunlade, and
  E.~Zhang, ``Accelerated high-resolution photoacoustic tomography via
  compressed sensing,'' \emph{Phys. Med. Biol.}, vol.~61, no.~24, p. 8908,
  2016.

\bibitem{haltmeier2016compressed}
M.~Haltmeier, T.~Berer, S.~Moon, and P.~Burgholzer, ``Compressed sensing and
  sparsity in photoacoustic tomography,'' \emph{J. Opt.}, vol.~18, no.~11, pp.
  114\,004--12pp, 2016.

\bibitem{haltmeier2017new}
M.~Haltmeier, M.~Sandbichler, T.~Berer, J.~Bauer-Marschallinger, P.~Burgholzer,
  and L.~Nguyen, ``A new sparsification and reconstruction strategy for
  compressed sensing photoacoustic tomography,'' \emph{arXiv:1801.00117}, 2017.

\bibitem{sandbichler2015novel}
M.~Sandbichler, F.~Krahmer, T.~Berer, P.~Burgholzer, and M.~Haltmeier, ``A
  novel compressed sensing scheme for photoacoustic tomography,'' \emph{SIAM J.
  Appl. Math.}, vol.~75, no.~6, pp. 2475--2494, 2015.

\bibitem{provost2009application}
J.~Provost and F.~Lesage, ``The application of compressed sensing for
  photo-acoustic tomography,'' \emph{IEEE Trans. Med. Imag.}, vol.~28, no.~4,
  pp. 585--594, 2009.

\bibitem{paltauf2007experimental}
G.~Paltauf, R.~Nuster, M.~Haltmeier, and P.~Burgholzer, ``Experimental
  evaluation of reconstruction algorithms for limited view photoacoustic
  tomography with line detectors,'' \emph{Inverse Probl.}, vol.~23, no.~6, pp.
  S81--S94, 2007.

\bibitem{paltauf2009weight}
G.~Paltauf, R.~Nuster, and P.~Burgholzer, ``Weight factors for limited angle
  photoacoustic tomography,'' \emph{Phys. Med. Biol.}, vol.~54, no.~11, p.
  3303, 2009.

\bibitem{antholzer2017deep}
S.~Antholzer, M.~Haltmeier, and J.~Schwab, ``Deep learning for photoacoustic
  tomography from sparse data,'' \emph{arXiv:1704.04587}, 2017, to appear in
  Inverse Probl Sci Eng.

\bibitem{pelt2014improving}
D.~M. Pelt and K.~J. Batenburg, ``Improving filtered backprojection
  reconstruction by data-dependent filtering,'' \emph{{IEEE} Trans. Image
  Process.}, vol.~23, no.~11, pp. 4750--4762, 2014.

\bibitem{wang2016perspective}
G.~Wang, ``A perspective on deep imaging,'' \emph{IEEE Access}, vol.~4, pp.
  8914--8924, 2016.

\bibitem{wang2016accelerating}
S.~Wang, Z.~Su, L.~Ying, X.~Peng, S.~Zhu, F.~Liang, D.~Feng, and D.~Liang,
  ``Accelerating magnetic resonance imaging via deep learning,'' in \emph{IEEE
  13th International Symposium on Biomedical Imaging (ISBI)}, 2016, pp.
  514--517.

\bibitem{chen2017lowdose}
H.~Chen, Y.~Zhang, W.~Zhang, P.~Liao, K.~Li, J.~Zhou, and G.~Wang, ``Low-dose
  {CT} via convolutional neural network,'' \emph{Biomed. Opt. Express}, vol.~8,
  no.~2, pp. 679--694, 2017.

\bibitem{jin2017deep}
K.~H. Jin, M.~T. McCann, E.~Froustey, and M.~Unser, ``Deep convolutional neural
  network for inverse problems in imaging,'' \emph{IEEE Trans. Image Process.},
  vol.~26, no.~9, pp. 4509--4522, 2017.

\bibitem{han2016deep}
Y.~Han, J.~J. Yoo, and J.~C. Ye, ``Deep residual learning for compressed
  sensing {CT} reconstruction via persistent homology analysis,''
  \emph{arXiv:1611.06391 [cs.CV]}, 2016.

\bibitem{wurfl2016deep}
T.~W{\"u}rfl, F.~C. Ghesu, V.~Christlein, and A.~Maier, ``Deep learning
  computed tomography,'' in \emph{International Conference on Medical Image
  Computing and Computer-Assisted Intervention}.\hskip 1em plus 0.5em minus
  0.4em\relax Springer, 2016, pp. 432--440.

\bibitem{zhang2016image}
H.~Zhang, L.~Li, K.~Qiao, L.~Wang, B.~Yan, L.~Li, and G.~Hu, ``Image prediction
  for limited-angle tomography via deep learning with convolutional neural
  network,'' arXiv:1607.08707, 2016.

\bibitem{rivenson2017deepb}
Y.~Rivenson, Z.~G{\"o}r{\"o}cs, H.~G{\"u}naydin, Y.~Zhang, H.~Wang, and
  A.~Ozcan, ``Deep learning microscopy,'' \emph{Optica}, vol.~4, no.~11, pp.
  1437--1443, 2017.

\bibitem{antholzer2018photoacoustic}
S.~Antholzer, M.~Haltmeier, R.~Nuster, and J.~Schwab, ``Photoacoustic image
  reconstruction via deep learning,'' in \emph{Photons Plus Ultrasound: Imaging
  and Sensing 2018}, vol. 10494.\hskip 1em plus 0.5em minus 0.4em\relax
  International Society for Optics and Photonics, 2018, p. 104944U.

\bibitem{allman2018photoacoustic}
D.~Allman, A.~Reiter, and M.~A.~L. Bell, ``Photoacoustic source detection and
  reflection artifact removal enabled by deep learning,'' \emph{IEEE Trans.
  Med. Imaging}, vol.~37, no.~6, pp. 1464--1477, 2018.

\bibitem{schwab2018deep}
J.~Schwab, S.~Antholzer, and M.~Haltmeier, ``Deep null space learning for
  inverse problems: Convergence analysis and rates,'' \emph{arXiv:1806.06137},
  2018.

\bibitem{adler2017solving}
J.~Adler and O.~{\"O}ktem, ``Solving ill-posed inverse problems using iterative
  deep neural networks,'' \emph{arXiv:1704.04058}, 2017.

\bibitem{hauptmann2018model}
A.~Hauptmann, F.~Lucka, M.~Betcke, N.~Huynh, J.~Adler, B.~Cox, P.~Beard,
  S.~Ourselin, and S.~Arridge, ``Model-based learning for accelerated,
  limited-view 3-d photoacoustic tomography,'' \emph{IEEE Trans. Med. Imaging},
  vol.~37, no.~6, pp. 1382--1393, 2018.

\bibitem{kelly2017deep}
B.~Kelly, T.~P. Matthews, and M.~A. Anastasio, ``Deep learning-guided image
  reconstruction from incomplete data,'' \emph{arXiv:1709.00584}, 2017.

\bibitem{gupta2018cnn}
H.~Gupta, K.~H. Jin, H.~Q. Nguyen, M.~T. McCann, and M.~Unser, ``Cnn-based
  projected gradient descent for consistent ct image reconstruction,''
  \emph{IEEE Trans. Med. Imag.}, vol.~37, no.~6, pp. 1440--1453, 2018.

\bibitem{kobler2017variational}
E.~Kobler, T.~Klatzer, K.~Hammernik, and T.~Pock, ``Variational networks:
  connecting variational methods and deep learning,'' in \emph{German
  Conference on Pattern Recognition}.\hskip 1em plus 0.5em minus 0.4em\relax
  Springer, 2017, pp. 281--293.

\bibitem{waibel2018reconstruction}
D.~Waibel, J.~Gr{\"o}hl, F.~Isensee, T.~Kirchner, K.~Maier-Hein, and
  L.~Maier-Hein, ``Reconstruction of initial pressure from limited view
  photoacoustic images using deep learning,'' in \emph{Photons Plus Ultrasound:
  Imaging and Sensing 2018}, vol. 10494.\hskip 1em plus 0.5em minus 0.4em\relax
  International Society for Optics and Photonics, 2018, p. 104942S.

\bibitem{XuWan03}
M.~Xu and L.~V. Wang, ``Analytic explanation of spatial resolution related to
  bandwidth and detector aperture size in thermoacoustic or photoacoustic
  reconstruction,'' \emph{Phys. Rev. E}, vol.~67, no.~5, pp.
  0\,566\,051--05\,660\,515 (electronic), 2003.

\bibitem{HalZan10}
M.~Haltmeier and G.~Zangerl, ``Spatial resolution in photoacoustic tomography:
  Effects of detector size and detector bandwidth,'' \emph{Inverse Probl.},
  vol.~26, no.~12, p. 125002, 2010.

\bibitem{haltmeier2016sampling}
M.~Haltmeier, ``Sampling conditions for the circular radon transform,''
  \emph{{IEEE} Trans. Image Process.}, vol.~25, no.~6, pp. 2910--2919, 2016.

\bibitem{burgholzer2007exact}
P.~Burgholzer, G.~J. Matt, M.~Haltmeier, and G.~Paltauf, ``Exact and
  approximate imaging methods for photoacoustic tomography using an arbitrary
  detection surface,'' \emph{Phys. Rev. E}, vol.~75, no.~4, p. 046706, 2007.

\bibitem{Treeby10}
B.~E. Treeby and B.~T. Cox, ``k-wave: Matlab toolbox for the simulation and
  reconstruction of photoacoustic wave-fields,'' \emph{J. Biomed. Opt.},
  vol.~15, p. 021314, 2010.

\bibitem{HriKucNgu08}
Y.~Hristova, P.~Kuchment, and L.~Nguyen, ``Reconstruction and time reversal in
  thermoacoustic tomography in acoustically homogeneous and inhomogeneous
  media,'' \emph{Inverse Probl.}, vol.~24, no.~5, p. 055006 (25pp), 2008.

\bibitem{FinHalRak07}
D.~Finch, M.~Haltmeier, and Rakesh, ``Inversion of spherical means and the wave
  equation in even dimensions,'' \emph{SIAM J. Appl. Math.}, vol.~68, no.~2,
  pp. 392--412, 2007.

\bibitem{Hal14}
M.~Haltmeier, ``Universal inversion formulas for recovering a function from
  spherical means,'' \emph{SIAM J. Math. Anal.}, vol.~46, no.~1, pp. 214--232,
  2014.

\bibitem{Hal13a}
------, ``Inversion of circular means and the wave equation on convex planar
  domains,'' \emph{Comput. Math. Appl.}, vol.~65, no.~7, pp. 1025--1036, 2013.

\bibitem{Kun07a}
L.~A. Kunyansky, ``Explicit inversion formulae for the spherical mean {R}adon
  transform,'' \emph{Inverse Probl.}, vol.~23, no.~1, pp. 373--383, 2007.

\bibitem{kuchment2008mathematics}
P.~Kuchment and L.~Kunyansky, ``Mathematics of thermoacoustic tomography,''
  \emph{Eur. J. Appl. Math.}, vol.~19, no.~2, pp. 191--224, 2008.

\bibitem{rosenthal2013acoustic}
A.~Rosenthal, V.~Ntziachristos, and D.~Razansky, ``Acoustic inversion in
  optoacoustic tomography: A review,'' \emph{Curr. Med. Imaging Rev.}, vol.~9,
  no.~4, pp. 318--336, 2013.

\bibitem{xu2005universal}
M.~Xu and L.~V. Wang, ``Universal back-projection algorithm for photoacoustic
  computed tomography,'' \emph{Phys. Rev. E}, vol.~71, no.~1, p. 016706, 2005.

\bibitem{burgholzer2007temporal}
P.~Burgholzer, J.~Bauer-Marschallinger, H.~Gr{\"u}n, M.~Haltmeier, and
  G.~Paltauf, ``Temporal back-projection algorithms for photoacoustic
  tomography with integrating line detectors,'' \emph{Inverse Probl.}, vol.~23,
  no.~6, pp. S65--S80, 2007.

\bibitem{HalPer15b}
M.~Haltmeier and S.~Pereverzyev, Jr., ``The universal back-projection formula
  for spherical means and the wave equation on certain quadric hypersurfaces,''
  \emph{J. Math. Anal. Appl.}, vol. 429, no.~1, pp. 366--382, 2015.

\bibitem{Nat12}
F.~Natterer, ``Photo-acoustic inversion in convex domains,'' \emph{Inverse
  Probl. Imaging}, vol.~6, no.~2, pp. 315--320, 2012.

\bibitem{deanben2012accurate}
X.~L. Dean-Ben, A.~Buehler, V.~Ntziachristos, and D.~Razansky, ``Accurate
  model-based reconstruction algorithm for three-dimensional optoacoustic
  tomography,'' \emph{{IEEE} Trans. Med. Imag.}, vol.~31, no.~10, pp.
  1922--1928, 2012.

\bibitem{haltmeier2017iterative}
M.~Haltmeier and L.~V. Nguyen, ``Analysis of iterative methods in photoacoustic
  tomography with variable sound speed,'' \emph{SIAM J. Imaging Sci.}, vol.~10,
  no.~2, pp. 751--781, 2017.

\bibitem{paltauf2002iterative}
G.~Paltauf, J.~A. Viator, S.~A. Prahl, and S.~L. Jacques, ``Iterative
  reconstruction algorithm for optoacoustic imaging,'' \emph{J. Opt. Soc. Am.},
  vol. 112, no.~4, pp. 1536--1544, 2002.

\bibitem{wang2012investigation}
K.~Wang, R.~Su, A.~A. Oraevsky, and M.~A. Anastasio, ``Investigation of
  iterative image reconstruction in three-dimensional optoacoustic
  tomography,'' \emph{Phys. Med. Biol.}, vol.~57, no.~17, p. 5399, 2012.

\bibitem{nguyen2014reconstruction}
L.~V. Nguyen, ``On a reconstruction formula for spherical radon transform: a
  microlocal analytic point of view,'' \emph{Anal. Math. Phys.}, vol.~4, no.~3,
  pp. 199--220, 2014.

\bibitem{SteUhl09}
P.~Stefanov and G.~Uhlmann, ``Thermoacoustic tomography with variable sound
  speed,'' \emph{Inverse Probl.}, vol.~25, no.~7, pp. 075\,011, 16, 2009.

\bibitem{haltmeier2010spatial}
M.~Haltmeier and G.~Zangerl, ``Spatial resolution in photoacoustic tomography:
  effects of detector size and detector bandwidth,'' \emph{Inverse Probl.},
  vol.~26, no.~12, p. 125002, 2010.

\bibitem{goodfellow2016deep}
I.~Goodfellow, Y.~Bengio, and A.~Courville, \emph{Deep Learning}.\hskip 1em
  plus 0.5em minus 0.4em\relax MIT Press, 2016.

\bibitem{han2017deep}
Y.~Han, J.~Gu, and J.~C. Ye, ``Deep learning interior tomography for
  region-of-interest reconstruction,'' \emph{arXiv:1712.10248}, 2017.

\bibitem{ronneberger2015unet}
O.~Ronneberge, P.~Fischer, and T.~Brox, ``U-net: Convolutional networks for
  biomedical image segmentation,'' in \emph{International Conference on Medical
  Image Computing and Computer-Assisted Intervention}, 2015, pp. 234--241.

\bibitem{ye2018deep}
J.~C. Ye, Y.~Han, and E.~Cha, ``Deep convolutional framelets: A general deep
  learning framework for inverse problems,'' \emph{SIAM J. Imaging Sci.},
  vol.~11, no.~2, pp. 991--1048, 2018.

\bibitem{mao2016image}
X.~Mao, C.~Shen, and Y.~Yang, ``Image restoration using very deep convolutional
  encoder-decoder networks with symmetric skip connections,'' in \emph{Advances
  in neural information processing systems}, 2016, pp. 2802--2810.

\bibitem{mangasarian1994backpropagation}
O.~L. Mangasarian and M.~V. Solodov, ``Backpropagation convergence via
  deterministic nonmonotone perturbed minimization,'' in \emph{Advances in
  Neural Information Processing Systems}, 1994, pp. 383--390.

\bibitem{tseng1998incremental}
P.~Tseng, ``An incremental gradient (-projection) method with momentum term and
  adaptive stepsize rule,'' \emph{SIAM J. Optim.}, vol.~8, no.~2, pp. 506--531,
  1998.

\bibitem{sidky2012convex}
E.~Y. Sidky, J.~H. J{\o}rgensen, and X.~Pan, ``Convex optimization problem
  prototyping for image reconstruction in computed tomography with the
  chambolle--pock algorithm,'' \emph{Phys. Med. Biol.}, vol.~57, no.~10, p.
  3065, 2012.

\bibitem{chambolle2011first}
A.~Chambolle and T.~Pock, ``A first-order primal-dual algorithm for convex
  problems with applications to imaging,'' \emph{J. Math. Imaging Vision},
  vol.~40, no.~1, pp. 120--145, 2011.

\bibitem{boink2018framework}
Y.~E. Boink, M.~J. Lagerwerf, W.~Steenbergen, S.~A. van Gils, S.~Manohar, and
  C.~Brune, ``A framework for directional and higher-order reconstruction in
  photoacoustic tomography,'' \emph{Phys. Med. \& Biol.}, vol.~63, no.~4, p.
  045018, 2018.

\bibitem{nguyen2018reconstruction}
L.~V. Nguyen and M.~Haltmeier, ``Reconstruction algorithms for photoacoustic
  tomography in heterogenous damping media,'' 2018, preprint.

\bibitem{chollet2015keras}
F.~Chollet, ``Keras,'' \url{https://github.com/fchollet/keras}, 2015, gitHub.

\bibitem{wang2004image}
Z.~Wang, A.~C. Bovik, H.~R. S., and E.~P. Simoncelli, ``Image quality
  assessment: from error visibility to structural similarity,'' \emph{IEEE
  Trans Image Process}, vol.~13, no.~4, pp. 600--612, 2004.

\bibitem{li2018nett}
H.~Li, J.~Schwab, S.~Antholzer, and M.~Haltmeier, ``{NETT}: Solving inverse
  problems with deep neural networks,'' \emph{arXiv:1803.00092}, 2018.

\end{thebibliography}
\end{document}